\documentclass[a4paper,12pt]{amsart}
\usepackage{amsfonts}
\usepackage{amsthm}
\usepackage{amssymb}
\usepackage{amsmath}
\usepackage{enumerate}
\usepackage{url}
\usepackage{color}
\usepackage{comment}
\usepackage{a4wide}
\usepackage{graphicx}
\usepackage{bbm}

\newtheorem{lemma}{Lemma}
\newtheorem{thm}{Theorem}
\newtheorem{defn}{Definition}

\newtheorem{thma}{Theorem}

\theoremstyle{remark}

\newcommand{\ve}{\varepsilon}
\newcommand {\e} {\varepsilon}

\allowdisplaybreaks

\begin{document}
\title{Lacunary sequences in analysis, probability and number theory}
\author[C.\ Aistleitner, I.\ Berkes and R.\ Tichy]{Christoph Aistleitner, Istv\'an Berkes and Robert Tichy}

\newcommand{\mods}[1]{\,(\mathrm{mod}\,{#1})}

\begin{abstract}
In this paper we present the theory of lacunary trigonometric sums and lacunary sums of dilated functions, from the origins of the subject up to recent develop\-ments. We describe the connections with mathematical topics such as equidistribution and discrepancy, metric number theory, normality, pseudorandomness, Diophantine equations, and the subsequence principle. In the final section of the paper we prove new results which provide necessary and sufficient conditions for the central limit theorem for subsequences, in the spirit of Nikishin's resonance theorem for convergence systems. More precisely, we characterize those sequences of random variables which allow to extract a subsequence satisfying a strong form of the central limit theorem.
\end{abstract}

\setcounter{tocdepth}{1}
\maketitle
\tableofcontents

\section{Introduction}

The word ``lacunary'' has its origin in the Latin \emph{lacuna} (ditch, gap), which is a diminutive form of \emph{lacus} (lake). Accordingly, a lacunary sequence is a sequence with gaps, and a lacunary trigonometric sum is a sum of trigonometric functions with gaps between the frequencies of consecutive summands. The origin of the theory of lacunary sums might lie in Weierstrass' famous example of a continuous, nowhere differentiable function (1872). Since then the subject has evolved into many very different directions, reflecting for example the emergence of modern measure theory and axiomatic probability theory in the early twentieth century, profound developments in harmonic analysis and Diophantine approximation, the establishment of ergodic theory as one of the key instruments of number theory, or the interest in notions of pseudo-randomness which are associated with the evolution of theoretical computer science. Throughout this paper we will be concerned with convergence/divergence properties of infinite trigonometric series
$$
\sum_{k=1}^\infty c_k \cos (2 \pi n_k x) \qquad \text{or} \qquad \sum_{k=1}^\infty c_k \sin (2 \pi n_k x),
$$
as well as with the asymptotic order and the distributional behavior of finite trigonometric sums
$$
\sum_{k=1}^N c_k \cos (2 \pi n_k x) \qquad \text{or} \qquad \sum_{k=1}^N c_k \sin (2 \pi n_k x)
$$
(the latter often in the simple case where $c_k \equiv 1$), and with their generalizations
$$
\sum_{k=1}^\infty c_k f (n_k x) \qquad \text{and} \qquad \sum_{k=1}^N c_k f (n_k x).
$$
Here $(c_k)_{k \geq 1}$ is a sequence of coefficients, and $(n_k)_{k \geq 1}$ is a sequence of positive integers (typically increasing), which satisfies some gap property such as the classical Hadamard gap condition
$$
\frac{n_{k+1}}{n_k} > q > 1, \qquad k \geq 1,
$$
or the ``large gap condition'' (also called ``super-lacunarity property'')
$$
\frac{n_{k+1}}{n_k} \to \infty, \qquad k \to \infty.
$$
Furthermore, $f$ is a 1-periodic function which is usually assumed to satisfy some regularity properties (such as being of bounded variation, being Lipschitz-continuous, etc.), and which for simplicity is usually assumed to be centered such that $\int_0^1 f(x)~dx=0$.\\

Early appearances of such lacunary sums include the following.\\

\begin{itemize}
 \item Sums of the form $\sum_{k=1}^N f(2^k x)$, where $f$ is an indicator function of a dyadic sub-interval of $[0,1]$, extended periodically with period 1. Borel used such sums in 1909 to show that almost all reals are ``normal''; more on this topic is contained in Section~\ref{sec:norm} below.\\
 \item Uniform distribution of sequences $(\{n_k x\})_{k \geq 1}$ in Weyl's seminar paper of 1916; more on this in Section~\ref{sec:UD}. Here and throughout the paper, we write $\{\cdot\}$ for the fractional part function.\\
 \item Kolmogorov's theorem on the almost everywhere convergence of lacunary trigonometric series if the sequence of coefficients is square-summable (1924), a result related to his Three Series Theorem for the almost sure convergence of series of independent random variables. Later it turned out that Kolmogorov's convergence theorem for trigonometric series actually remains true without any gap condition whatsoever, a result which was widely believed to be ``too good to be true'' before being established by Carleson~\cite{carle} in 1966. More on this in Section~\ref{sec:UD}.\\
 \item Foundational work on the distribution of normalized lacunary trigonometric sums, in particular the central limit theorems of Kac (1946) and Salem and Zygmund (1947), and the laws of the iterated logarithms of Salem and Zygmund (1950) and of Erd\H os and G\'al (1955). More on this in Sections~\ref{sec:CLT} and~\ref{sec:LIL}.\\
 \end{itemize}

A fundamental observation is that the unit interval, equipped with Borel sets and Lebesgue measure, forms a probability space, and that consequently a sequence of functions such as $(\cos (2 \pi n_k x))_{k \geq 1}$ or $(f(n_k x))_{k \geq 1}$ can be viewed as a sequence of random variables over this space; if $f$ is 1-periodic and if $(n_k)_{k \geq 1}$ is a sequence of positive integers then these random variables are identically distributed, but in general they are not independent. However, under appropriate circumstances the gap condition which is imposed upon $(n_k)_{k \geq 1}$ can ensure that these random variables have a low degree of stochastic dependence. Consequently lacunary sums often mimic the behavior of sums of independent and identically distributed random variables. This viewpoint was in particular taken by Steinhaus, Kac, and Salem and Zygmund in their fundamental work on the subject. In a particularly striking situation, the dyadic functions considered by Borel actually turn out to be a version of a sequence of Bernoulli random variables which are truly stochastically independent; accordingly, Borel's result on the normality of almost all reals is nowadays usually read as the historically very first version of the strong law of large numbers in probability theory.\\

When taking this probabilistic viewpoint, the theory of lacunary sums could be seen as a particular segment of the much wider field of the theory of weakly dependent random systems in probability theory, which is associated with notions such as mixing, martingales, and short-range dependence. However, it should be noted that the precise dependence structure in a lacunary function system $(f(n_k x))_{k \geq 1}$ is controlled by the analytic properties of the function $f$, in conjunction with arithmetic properties of the sequence $(n_k)_{k \geq 1}$. It is precisely this interplay between probabilistic, analytic and arithmetic aspects which makes the theory of lacunary sums so interesting, so challenging and so rewarding. In the following sections we want to illustrate some instances of these phenomena in more detail.

\section{Uniform distribution and discrepancy} \label{sec:UD}

Let $(x_n)_{n \geq 1}$ be a sequence of real numbers in the unit interval. Such a sequence is called uniformly distributed modulo one (in short: u.d.\ mod 1) if
\begin{equation} \label{ud}
\lim_{N \to \infty} \frac{1}{N} \sum_{n=1}^N \mathbbm{1}_A (x_n) =  \lambda(A)
\end{equation}
for all sub-intervals $A \subset [0,1]$ of the unit interval. The word ``equidistributed'' is also used for this property, synonymously with ``uniformly distributed modulo one''. In this definition, and in the sequel, $\mathbbm{1}$ denotes an indicator function, and $\lambda$ denotes Lebesgue measure. In informal language, this definition means that a sequence is u.d.\ mod 1 if every interval $A$ asymptotically receives its fair share of elements of the sequence, which is proportional to the length of the interval. Note that (for example as a consequence of the Glivenko--Cantelli theorem) for a sequence of independent, uniformly $(0,1)$-distributed random variables $(U_n)_{n \geq 1}$ one has
$$
\frac{1}{N} \sum_{n=1}^N \mathbbm{1}_A (U_n) =  \lambda(A) \qquad \text{almost surely}
$$
for all intervals $A \subset [0,1]$, so that in a vague sense uniform distribution of a deterministic sequence can be interpreted in the sense that the sequence shows ``random'' behavior; more on this aspect in Section~\ref{sec:norm} below. Uniform distribution theory can be said to originate with Kronecker's approximation theorem and with work of Bohl, Sierpi\'nski and Weyl on the sequence $(\{n \alpha\})_{n \geq 1}$ for irrational $\alpha$. However, the theory only came into its own with Hermann Weyl's~\cite{weyl} seminal paper of 1916. Among many other fundamental insights, Weyl realized that Definition \eqref{ud}, which in earlier work had only be read in terms of counting points in certain intervals, can be interpreted in a ``functional'' way and can equivalently be written as
\begin{equation} \label{ud2}
\lim_{N \to \infty} \frac{1}{N} \sum_{n=1}^N f (x_n) =  \int_0^1 f(x)~dx
\end{equation}
for all continuous functions $f$. This viewpoint suggests that uniformly distributed sequences can be used as quadrature points for numerical integration; in the multi-dimensional setting and together with quantitative error estimates this observation forms the foundation of the so-called Quasi-Monte Carlo integration method, a concept which today forms a cornerstone of numerical methods in quantitative finance and other fields of applied mathematics (more on this below). Furthermore, Weyl realized that the indicator functions in \eqref{ud} or the continuous functions in \eqref{ud2} could also be replaced by complex exponentials, as a consequence of the Weierstrass approximation theorem; thus by the famous Weyl Criterion a sequence is u.d.\ mod 1 if and only if
$$
\lim_{N \to \infty} \frac{1}{N} \sum_{n=1}^N  e^{2 \pi i h x_n} =  0
$$
for all fixed non-zero integers $h$, thereby tightly connecting uniform distribution theory with the theory of exponential sums.\\

For the particular sequence $(\{n \alpha\})_{n \geq 1}$ it can be easily seen from the Weyl criterion that this sequence is u.d.\ mod 1 if and only if $\alpha \not\in \mathbb{Q}$. However, for other parametric sequences of the form $(\{n_k \alpha\})_{k \geq 1}$ the situation is much more difficult, and in general it is completely impossible to determine whether for some particular value of $\alpha$ the sequence is u.d.\ or not. It turns out that in a \emph{metric} sense the situation is quite different. Metric number theory arose after the clarification of the concept of real numbers, the realization that the reals drastically outnumber the integers and the rationals, and the development of modern measure theory. Loosely speaking, the purpose of metric number theory is to determine properties which hold for a set of reals which is ``typical'' with respect to a certain measure; here ``typical'' means that the measure of the complement is small. In the present paper the measure under consideration will always be the Lebesgue measure, and a set of reals will be considered typical if its complement has vanishing Lebesgue measure; however, metric number theory has for example also been intensively studied with respect to the Hausdorff dimension or other fractal measures. \\

Returning to Weyl's results, what he proved in the metric setting is the following. For every sequence of distinct integers $(n_k)_{k \geq 1}$, the sequence $(\{n_k \alpha\})_{k \geq 1}$ is u.d.\ mod 1 for (Lebesgue-) almost all reals $\alpha$. In other words, even if we cannot specify the set of $\alpha$'s for which uniform distribution holds, at least we know that the set of such $\alpha$'s has full Lebesgue measure. It is amusing that after formulating the result, Weyl continues to write:

\vskip3mm
\begin{quote}
 \emph{Wenn ich nun freilich glaube, da{\ss} man den Wert solcher S\"atze, in denen eine unbestimmte Ausnahmemenge vom Ma{\ss}e 0 auftritt, nicht eben hoch einsch\"atzen darf, m\"ochte ich diese Behauptung hier doch kurz begr\"unden.}\footnote{Even if I think that the value of theorems, which contain an unspecified exceptional set of measure zero, is not particularly high, I still want to give a short justification.}
\end{quote}
\vskip3mm

One should recall that Weyl's paper was written in a time of intense conflict of formalists vs.\ constructivists (with Weyl favoring the latter ones), and only very briefly after the notion of a set of zero (Lebesgue) measure had been introduced at all. Today, Weyl's theorem is seen as one of the foundational results of metric number theory, together with the work of Borel, Koksma, Khinchin and others.\\

While uniform distribution modulo one is a qualitative asymptotic property, it is natural that one is also interested in having a corresponding quantitative concept which applies to finite sequences (or finite truncations of infinite sequences). Such a concept is the discrepancy of a sequence, which is defined by
$$
D_N(x_1, \dots, x_N) = \sup_{A \subset [0,1]} \left|\frac{1}{N} \sum_{n=1}^N \mathbbm{1}_A (x_n) - \lambda(A) \right|.
$$
Here the supremum is taken over all sub-intervals $A \subset [0,1]$, and it is easy to see that an infinite sequence $(x_n)_{n \geq 1}$ is u.d.\ mod 1 if and only if the discrepancy $D_N(x_1, \dots,x_N)$ tends to 0 as $N \to \infty$. With a slight abuse of notation, we will write throughout the paper $D_N(x_n)=D_N(x_1, \dots, x_N)$ for the discrepancy of the first $N$ elements of an infinite sequence $(x_n)_{n \geq 1}$. From a probabilistic perspective, the discrepancy is a variant of the (two-sided) Kolmogorov--Smirnov statistic, where one tests the empirical distribution of the point set $x_1, \dots, x_N$ against the uniform distribution on $[0,1]$. Without going into details, we note that $D_N(x_1, \dots, x_N)$ can be bounded above in terms of exponential sums by the Erd\H os--Tur\'an inequality, and that the error when using $x_1, \dots, x_N$ as a set of quadrature points to approximate $\int_0^1 f(x)~dx$ by $\frac{1}{N} \sum_{n=1}^N f(x_n)$ can be bounded above by Koksma's inequality in terms of the variation of $f$ and the discrepancy $D_N$; for details see the monographs~\cite{dts,kn}, which contain all the basic information on uniform distribution theory and discrepancy. See also~\cite{mont} for a discussion of equidistribution and discrepancy from the viewpoint of analytic number theory, and~\cite{lemi,leop,novw} for expositions which put particular emphasis on the numerical analysis aspects. \\

Weyl's metric result from above can be written as
$$
\lim_{N \to \infty} D_N(\{n_k \alpha\}) = 0 \qquad  \text{for almost all $\alpha$},
$$
for any sequence $(n_k)_{k \geq 1}$ of distinct itegers. Strikingly, the precise answer to the corresponding quantitative problem is still open more than a hundred years later. It is known that for every strictly increasing sequence of integers $(n_k)_{k \geq 1}$ one has
\begin{equation} \label{discr_upper}
D_N(\{n_k \alpha\}) = O \left( \frac{(\log N)^{3/2+\ve}}{\sqrt{N}} \right) \qquad \text{for almost all $\alpha$}.
\end{equation}
This is a result of R.C.\ Baker~\cite{baker}, who improved earlier results of Cassels~\cite{cass} and of Erd\H os and Koksma~\cite{ek} by using Carleson's celebrated convergence theorem in the form of the Carleson--Hunt inequality~\cite{hunt}. In his paper Baker wrote that

\vskip3mm
\begin{quote}
 [\dots] probably the exponent $3/2 + \varepsilon$ could be replaced by $\ve$ [\dots]
\end{quote}
\vskip3mm

but it turned out that this is not actually the case. Instead, Berkes and Philipp~\cite{bpt} constructed an example of an increasing integer sequence $(n_k)_{k \geq 1}$ for which
\begin{equation} \label{berkphil}
\limsup_{N \to \infty} \frac{\left| \sum_{k=1}^N \cos (2 \pi n_k x) \right|}{\sqrt{N \log N}} = +\infty \qquad \text{for almost every $x$}.
\end{equation}
By the Erd\H os--Tur\'an inequality this gives a corresponding lower bound for the discrepancy, which implies that the optimal exponent of the logarithmic term in an upper bound of the form \eqref{discr_upper} has to be at least 1/2. But the actual size of this optimal exponent, one of the most fundamental problems in metric discrepancy theory, still remains open. Note that for pure cosine-sums $\sum_{k=1}^N \cos (2 \pi n_k x)$ it is easily seen that one has a metric upper bound with exponent $1/2 + \varepsilon$ in the logarithmic term; this follows from the orthogonality of the trigonometric system, together with Carleson's inequality and the Chebyshev inequality. Thus, in connection with \eqref{berkphil}, the optimal upper bound in a metric estimate for pure cosine sums is known. For sums $\sum_{k=1}^N f (n_k x)$ with $f$ being a 1-periodic function of bounded variation, the optimal exponent also is $1/2 + \varepsilon$, but this is a much deeper result than the one for the pure cosine case, and was established only recently in~\cite{abs,lr}. By Koksma's inequality, an upper bound for the discrepancy implies an upper bound for sums of function values for a (fixed) function of bounded variation, but the opposite is not true. So while the case of a fixed function $f$ is solved and is an important test case for the discrepancy, the problem of the discrepancy itself (which requires a supremum over a whole class of test functions) is more involved and remains open.

\section{Arithmetic effects: Diophantine equations and sums of common divisors}\label{sec:arithm}

One of the most classical tools of probability theory is the calculation of expectations, variances, and higher moments of sums of random variables. Due to trigonometric identities such as
\begin{equation} \label{trig_id}
\cos a \cos b = \frac{\cos(a+b) + \cos(a-b)}{2},
\end{equation}
the calculation of moments of sums of trigonometric functions (with integer frequencies) reduces to a counting of solutions of certain Diophantine equations. Indeed, while the first and second moments
$$
\int_0^1 \sum_{k=1}^N \cos(2 \pi n_k x)~dx = 0 \qquad \text{and} \qquad \int_0^1 \left(\sum_{k=1}^N \cos(2 \pi n_k x)\right)^2 dx = \frac{N}{2}
$$
are trivial and do not depend on the particular sequence $(n_k)_{k \geq 1}$ (as long as the elements of the sequence are assumed to be distinct), interesting arithmetic effects come into play when one has to compute higher moments, and it can be clearly seen how the presence of a gap condition leads to a behavior of the moments which is similar to that of sums of independent random variables. More precisely, assume that we try to calculate
$$
\int_0^1 \left(\sum_{k=1}^N \cos(2 \pi n_k x)\right)^m dx
$$
for some integer $m \geq 3$. By \eqref{trig_id} this can be written as a sum
$$
2^{-m} ~\sum_{\pm}~ \sum_{1 \leq k_1, \dots, k_m \leq N} \mathbbm{1} \left(\pm n_{k_1} \pm \dots \pm n_{k_m} =0 \right).
$$
Here the first sum is meant as a sum over all positive combinations of ``+'' and ``-'' signs inside the indicator function at the end. Now assume that, for simplicity, we consider the particular sequence $n_k = 2^k,~k \geq 1$, which is a prototypical example of a sequence satisfying the Hadamard gap condition. Then the majority of the solutions of
$$
\pm n_{k_1} \pm \dots \pm n_{k_m} =0
$$
arise from cancellation in a pairwise way; that is, after a suitable re-ordering of the indices, one has $k_1 = k_2, ~k_3 = k_4, \dots, ~k_{m-1} = k_m,$ and suitable opposing signs.  There are further solutions, coming for example from the fact that $n_{k+1} - n_k - n_k=0$ for all $k$, but one can show that they play an overall insignificant role. Thus one is led to calculating the combinatorial quantity which stems from the solutions with pairwise cancellation, which turns out to be exactly the same combinatorial quantity that arises when calculating an $m$-th moment of a sum of independent random variables. Thus the moments of the trigonometric lacunary sum converge to those of a suitable Gaussian distribution, which gives rise to the classical limit theorems for lacunary trigonometric sums. The situation is more delicate if one only has the Hadamard gap condition $n_{k+1}/n_k > q > 1$ rather than exact exponential growth, and again more delicate if one considers a sum of dilated functions $\sum f(n_k x)$ instead of a pure trigonometric sum, but the principle described here is very powerful also in these more general situations. For a long time this was the key ingredient in most of the proofs of limit theorems for lacunary sums; see for example~\cite{egot,kac,plt,sz1947,taka,weiss}. A different method is based on the approximation of a lacunary sum by a martingale difference; here the ``almost independent'' behavior is not captured by controlling the moments of the sum, but in the fact that later terms of the sum (functions with high frequency) oscillate quickly in small regions where earlier summands (functions with much lower frequency) are essentially constant. As far as we can say, this method was first used in the context of lacunary sums by Berkes~\cite{berk1} and, independently, by Philipp and Stout~\cite{psa}. We will come back to this topic in Section~\ref{sec:CLT}. \\

Broadly speaking, the ``almost independent'' behavior of sums of dilated functions breaks down when the lacunarity condition is relaxed. Many papers have been devoted to this effect; see in particular~\cite{berk_an, berk_cl, erd_o, murai}. In order to maintain the ``almost independent'' behavior of the sum, there are two natural routes to take. On the one hand, one could randomize the construction of the sequence $(n_k)_{k \geq 1}$, and assume that the undesired effects disappear almost surely with respect to the underlying probability measure -- it turns out that this is a very powerful method, and we will come back to it in Section~\ref{sec:rand} below. On the other hand, when adapting the viewpoint that the ``almost independence'' property is expressed in the small number of solutions of certain Diophantine equations, one could try to compensate the weaker growth assumption by stronger arithmetic assumptions. A prominent example of a class of sequences for which the latter approach has been very successfully used are the so-called Hardy--Littlewood--P\'olya sequences, which consist of all the elements of the multiplicative semigroup generated by a finite set of primes, sorted in increasing order. These sequences are in several ways a natural analogue of lacunary sequences; note that the sequence $(2^k)_{k \geq 1}$ actually also falls into this framework by consisting of all elements of the semigroup generated by a single prime. Such sequences generated by a finite set of primes have attracted the attention of number theorists again and again, a particularly interesting  instance being F\"urstenberg's~\cite{furst} paper on disjointness in ergodic theory. It is known that Hardy--Littlewood--P\'olya sequences (if generated by two or more primes) grow sub-exponentially, and the precise (only slightly sub-exponential) growth rates are known (Tijdeman~\cite{tijd}). What is more striking (and a much deeper fact) is that also the number of solutions of the relevant linear Diophantine equations can be bounded efficiently -- this is Schmidt's celebrated Subspace Theorem~\cite{schmidt} in a quantitative form such as that of Evertse, Schlickewei and Schmidt~\cite{ess} or Amoroso and Viada~\cite{av}. By a combination of the (slightly weaker) growth condition with the (strong) arithmetic information available for Hardy--Littlewood--P\'olya sequences, much of the machinery that is used for Hadamard lacunary sequences can be rescued for this generalized setup; see~\cite{phil_bak} as well as~\cite{abt_perm,bpt_emp,fuk_nak}.\\

We briefly come back to the case of sums of dilated functions $\sum f(n_k x)$ without the presence of a growth condition on $(n_k)_{k \geq 1}$. We assume for simplicity that $\int_0^1 f(x)~dx = 0$, so trivially
$$
\int_0^1 \sum_{k=1}^N f(n_k x) ~dx = 0,
$$
but already the calculation of the variance
\begin{equation} \label{var_calc}
\int_0^1 \left( \sum_{k=1}^N f(n_k x) \right)^2~dx
\end{equation}
is in general quite non-trivial. If $f(x) = \cos (2 \pi x)$, then one can simply use the orthogonality of the trigonometric system. If $f$ is a more general function, then one can still express $f$ by its Fourier series, expand the square and integrate, and thus translate the problem of calculating \eqref{var_calc} into a problem of counting the solutions of certain linear Diophantine equations. When carrying out this approach, one is naturally led to the problem of estimating a certain sum involving greatest common divisors. For example, assume that $f(x) = \{x\} - 1/2$. In this case a classical formula (first stated by Franel and first proved by Landau) asserts that
$$
\int_0^1 f(mx) f(nx)~ dx = \frac{1}{12} \frac{(\gcd(m,n))^2}{m n},
$$
and consequently
$$
\int_0^1 \left( \sum_{k=1}^N (\{n_k x\}-1/2) \right)^2 dx =  \frac{1}{12}  \sum_{1 \leq k,\ell \leq N} \frac{(\gcd(n_k,n_\ell))^2}{n_k n_\ell}.
$$
The sum on the right-hand side of this equation is called a GCD sum. A similar identity holds for example for the Hurwitz zeta function $\zeta(1-\alpha,\cdot)$, where
$$
\int_0^1 \zeta(1-\alpha,\{mx\}) \zeta(1-\alpha,\{n x\}) ~dx = 2 \Gamma(\alpha)^2 \frac{\zeta(2\alpha)}{(2 \pi)^{2 \alpha}}  \frac{(\gcd(m,n))^{2\alpha}}{(m n)^{\alpha}}
$$
for $\alpha>1/2$, thus leading to a GCD sum
$$
\sum_{1 \leq k,\ell \leq N} \frac{(\gcd(n_k,n_\ell))^{2\alpha}}{(n_k n_\ell)^\alpha}
$$
with parameter $\alpha$. If $f(x)$ is a general 1-periodic function, then one usually does not obtain such a nice exact representation of the variance of a sum of dilated function values, but typically the variance \eqref{var_calc} can be bounded above by a GCD sum, which together with Chebyshev's inequality and the Borel--Cantelli lemma allows to make a statement on the almost everywhere asymptotic behavior of a sum of dilated function values. \\

This connection between sums of dilated functions and GCD sums is explained in great detail in Chapter 3 in Harman's monograph on \emph{Metric Number Theory}~\cite{harman}, where mainly the context of metric Diophantine approximation is treated (see also~\cite{gal,koks}). Recently this connection has also led to a solution of the problem of the almost everywhere convergence of series of dilated functions. Recall that Carleson's theorem~\cite{carle} asserts that the series
$$
\sum_{k=1}^\infty c_k \cos (2 \pi n_k x)
$$
is almost everywhere convergent provided that $\sum_k c_k^2 < \infty$. It is natural to ask which assumption on the sequence of coefficients $(c_k)_{k \geq 1}$ is necessary to ensure the almost everywhere convergence of the more general series
\begin{equation} \label{this_ser}
\sum_{k=1}^\infty c_k f (n_k x),
\end{equation}
under some regularity assumptions on $f$.
Gaposhkin~\cite{gapo2,gapo1} obtained some partial results, but a satisfactory understanding of the problem was only achieved very recently, when the connection with GCD sums was fully understood and optimal upper bounds for such sums were obtained. Exploiting this connection with GCD sums, it was shown in~\cite{abs,lr} that for 1-periodic $f$ which is of bounded variation on $[0,1]$ the series \eqref{this_ser} is almost everywhere convergent provided that
$$
\sum_{k=3}^\infty c_k^2 (\log \log k)^\gamma < \infty
$$
for some $\gamma > 2$, and this result is optimal in the sense that the same assumption with $\gamma=2$ would not be sufficient. In~\cite{absw} it was shown that for the class $C_\alpha$ of 1-periodic square integrable functions $f$ with Fourier coefficients $a_j, b_j$ satisfying
$$ a_j=O(j^{-\alpha}), \quad b_j=O(j^{-\alpha}) $$
for $1/2< \alpha<1$, a sharp criterion for the almost everywhere convergence of  \eqref{this_ser} is that \begin{equation}\label{aaa2}
\sum_{k=1}^\infty c_k^2 \exp \left(\frac{K(\log k)^{1-\alpha}}{(\log\log k)^\alpha}\right)<\infty
\end{equation}
with a suitable $K=K(\alpha)$.
In the case of 1-periodic Lipschitz $\alpha$ functions $f$, Gaposhkin~\cite{gapo1} proved that for $\alpha>1/2$, the series  \eqref{this_ser} converges a.e.\  under $\sum_k c_k^2 < \infty$ (just like in the case of Carleson's theorem) and Berkes~\cite{berkes98} showed that this result is sharp, i.e.\ for $\alpha=1/2$ the exact analogue of Carleson's theorem is not valid. No sharp convergence criteria exists in the case $0<\alpha \le 1/2$; for sufficient criteria for see Gaposhkin~\cite{gapo2}. See also~\cite{alip,bewe1,gap1966} for a general discussion and several further results for the convergence of series $\sum_k c_k f(n_kx)$. \\

For general periodic $f\in L^2$ the direct connection between the integral (\ref{var_calc}) and GCD sums breaks down,  but upper bounds for (\ref{var_calc}) as well as for
\begin{equation}\label{ckfnkx}
\int_0^1 \left( \sum_{k=1}^N c_k f(n_kx)\right)^2 \, dx
\end{equation}
can be given in terms of the coefficients $c_k$, of the Fourier coefficients of $f$, and arithmetic functions such as
$d(n)=\sum_{d|n} 1$, $\sigma_s (n)= \sum_{d|n} d^s$, or the Erd\H{o}s-Hooley function $\Delta(n)=\sup_{u\in {\mathbb R}}\sum_{d|n, u\le d\le eu} 1$.  See Koksma~\cite{koks51,koks53}, Weber~\cite{weber2011}, and Berkes and Weber~\cite{bewe2, bewe3}.  A typical example (see~\cite{weber2011}) is the bound
$$ \int_0^1 \left( \sum_{k\in H} c_k f(kx)\right)^2 dx \le \left(\sum_{\nu=1}^\infty a_\nu^2 \Delta(\nu)\right) \sum_{k\in H} c_k^2 d(k)$$
valid for any set $H$ of disjoint positive integers lying in some interval $[e^r, e^{r+1}]$, $r\ge 1$. Here $a_k$ are the complex Fourier coefficients of $f$. Using standard methods, such bounds lead easily to a.e.\ convergence criteria for sums $\sum_k c_k f(kx)$, see the papers cited above. \\

In Wintner~\cite{wi} it was proved that if $f$ is a periodic $L^2$ function with Fourier coefficients $a_k, b_k$, then the series $\sum_k c_k f(kx)$ converges in $L^2$ norm for all coefficient sequences $(c_k)_{k\ge 1}$ satisfying $\sum_k c_k^2<\infty$ if and only if the functions defined by the Dirichlet series
$$ \sum_{k=1}^\infty a_k k^{-s}, \qquad \sum_{k=1}^\infty b_k k^{-s},$$
are bounded and regular in the half plane $\Re (s)>0$.
There is also a remarkable connection between the maximal order of magnitude of GCD sums with the order of extreme values of the Riemann zeta function in the critical strip; see~\cite{bs,dlbt,hilb,sound}.\\

Naturally, estimating the integral (\ref{var_calc}) provides important information also on the asymptotic behavior of averages
\begin{equation}\label{aver1}
\frac{1}{N} \sum_{k=1}^N f(n_kx).
\end{equation}
By the Weyl equidistribution theorem, for any 1-periodic $f$ with bounded variation in $(0, 1)$ we have
\begin{equation} \label{aver2}
\lim_{N\to\infty} \frac{1}{N} \sum_{k=1}^N f(kx)=\int_0^1 f(x)\, dx \quad \text{a.e.}
\end{equation}
(actually for every irrational $x$).
Khinchin~\cite{kh1923} conjectured that (\ref{aver2}) holds for every  1-periodic Lebesgue integrable $f$  as well. This conjecture remained open for nearly 50 years and was finally disproved by Marstrand~\cite{mars}. An example for a periodic integrable $f$  and a sequence $(n_k)_{k \geq 1}$ of positive integers such that the averages (\ref{aver1}) do not converge almost everywhere had already been given earlier by Erd\H{o}s~\cite{erd_49}.
On the other hand, Koksma~\cite{koks53} proved that (\ref{aver2}) holds if $f\in L^2$ and the Fourier coefficients $a_k$, $b_k$ of $f$ satisfy
$$ \sum_{k=1}^\infty \left((a_k^2+b_k^2) \sum_{d|k} \frac{1}{d} \right) <\infty,$$
and Berkes and Weber~\cite{bewe3} proved that the last condition is optimal. No similarly sharp criteria are known in the case $f\in L^1$. For further results related to the Khinchin conjecture, see~\cite{bakk,beck,bewe3,buc,nair}.

\section{The central limit theorem for lacunary sequences} \label{sec:CLT}

Salem and Zygmund~\cite{sz1947} proved the first central limit theorem (CLT) for lacunary trigonometric sums. More specifically, they showed that for any integer sequence $(n_k)_{k \geq 1}$ satisfying the Hadamard gap condition one has
$$
\lim_{N \to \infty} \lambda \left( x \in (0,1):~\sum_{k=1}^N \cos (2 \pi n_k x) \leq t \sqrt{N/2} \right) = \Phi(t),
$$
where $\Phi$ denotes the standard normal distribution. Note that
$$
\int_0^1 \left( \sum_{k=1}^N \cos(2 \pi n_k x) \right)^2 dx = \frac{N}{2},
$$
so the result above contains the ``correct'' variance for the limit distribution, exactly as it should also be expected in the truly independent case. This result has been significantly strengthened since then; for example, Philipp and Stout~\cite{psa} showed that under the Hadamard gap condition the function
$$
S(t,x) = \sum_{k \leq t} \cos (2 \pi n_k x),
$$
considered as a stochastic process over the space $([0,1],\mathcal{B}[0,1],\lambda)$, is a small perturbation of a Wiener process, a characterization which allows to deduce many fine asymptotic results for this sum. It is also known that the central limit theorem for pure trigonometric lacunary sums remains valid under a slightly weaker gap condition than Hadamard's: as Erd\H os~\cite{erd_o} proved, it is sufficient to assume that $n_{k+1}/n_k \geq 1 + c k^{-\alpha}$, $\alpha<1/2$,  while such an assumption with $\alpha=1/2$ is not sufficient.\\

The whole situation becomes very different when the cosine-function is replaced by a more general 1-periodic function, even if it is such a well-behaved one as a trigonometric polynomial. For example, consider
\begin{equation} \label{triv}
f(x) = \cos (2 \pi x) - \cos(4 \pi x), \qquad n_k = 2^k,~k \geq 1.
\end{equation}
In this case the lacunary sum is telescoping, and it can be immediately seen that there cannot be a non-trivial limit distribution. A more delicate example is attributed to Erd\H os and Fortet\footnote{The Erd\H os--Fortet example is first mentioned in print in a paper of Salem and Zygmund~\cite{sz1948}. They mention the example without proof, and write: ``This remark is essentially due to Erd\H{o}s.''. Later the example was mentioned in a paper of Kac~\cite{kac49}, who wrote: ``It thus came as a surprise when simultaneously and independently of each other, Erd\H{o}s and Fortet constructed an example showing that the limit [\dots] need not be Gaussian'', with a footnote: ``In Salem and Zygmund this example is erroneously credited to Erd\H os alone.'' No proof is given in Kac's paper either, but he writes: ``Details will be given in [a forthcoming] paper by Erd\H{o}s, Ferrand, Fortet and Kac''. Such a joint paper never appeared.}, and goes as follows. Let
$$
f(x) = \cos (2 \pi x) + \cos (4 \pi x), \qquad n_k = 2^k-1,~k \geq 1.
$$
Then it can be shown that $N^{-1/2} \sum_{k=1}^N f(n_k x)$ does indeed have a limit distribution, but one which is actually non-Gaussian. More precisely, for this example one has
$$
\lim_{N \to \infty} \lambda \left( x \in (0,1):~\sum_{k=1}^N f (n_k x) \leq t \sqrt{N/2} \right) = \frac{1}{\sqrt{\pi}} \int_0^1 \int_{-\infty}^{t/2 |\cos(\pi s)|} e^{-u^2} du ds.
$$
Thus the limit distribution in this case is a so-called ``variance mixture Gaussian'', which can be seen as a normal distribution whose variance is a function rather than a constant. This limiting behavior can be explained from the observation that
\begin{equation} \label{kac1}
f(n_k x) = \cos ( (2^{k+1}-2) \pi x) + \cos ( (2^{k+2}-4) \pi x)
\end{equation}
and
\begin{equation} \label{kac2}
f(n_{k+1} x) = \cos ( (2^{k+2}-2) \pi x) + \cos ( (2^{k+3}-4) \pi x).
\end{equation}
Combining the second term on the right-hand side of \eqref{kac1} with the first term on the right-hand side of \eqref{kac2} we obtain
$$
\cos ( (2^{k+2}-4) \pi x) + \cos ( (2^{k+2}-2) \pi x) = 2 \cos(\pi x) \cos((2^{k+2}-3)\pi x),
$$
so the whole lacunary sum $\sum_{k=1}^N f(n_k x)$ can essentially be written as $2 \cos (\pi x)$ multiplied with a pure cosine lacunary sum. This is exactly what the ``variance mixture Gaussian'' indicates: the limit distribution is actually that of $2\cos (\pi x)$ independently multiplied with a Gaussian.
The failure of a Gaussian central limit theorem in the example above can be seen as a consequence of the fact that the Diophantine equation
$$
n_{k+1} - 2 n_k = 1
$$
possesses many solutions $k$ for this particular choice of sequence. Equipped with this observation, one could readily construct similar examples with other trigonometric polynomials $f$, and other variance mixture Gaussians as limit distributions, by creating situations where there are many solutions $k,\ell$ to
\begin{equation} \label{dio_equ}
a n_k - b n_\ell = c
\end{equation}
for some fixed $a,b,c$. However, interestingly, a special role is played by such equations when $c$ has the particular value $c=0$; very roughly speaking, solutions of the equation for $c=0$ effect only the limiting variance (in a Gaussian distribution), but not the structure of the limiting distribution itself. This is visible in a paper of Kac~\cite{kac}, who studied the sequence $n_k=2^k,~k \geq 1$, where indeed the only equations that have many solution are of the form $2^m n_k - n_\ell = 0$ for some $m$ (the solutions being $\ell  = k+m$). Kac proved that for this sequence and any 1-periodic $f$ of bounded variation and zero mean one has
\begin{equation} \label{kac1a}
\lim_{N \to \infty} \lambda \left( x \in (0,1):~\sum_{k=1}^N f (n_k x) \leq t \sigma_f \sqrt{N} \right) = \Phi(t)
\end{equation}
with a limiting variance $\sigma_f^2$, provided that
\begin{equation} \label{kac2a}
\sigma_f^2:= \int_0^1 f^2(x) ~dx + 2 \sum_{m=1}^\infty \int_0^1 f(x) f(2^m x)~dx \neq 0.
\end{equation}
Thus in this case the limit distribution is always a Gaussian, and the failure of the trivial example in \eqref{triv} to produce such a Gaussian limit comes from the fact that the limiting variance is degenerate.\\

These observations show that there is a delicate interplay between arithmetic, analytic and probabilistic effects; in particular, it is obviously not only the order of growth of $(n_k)_{k \geq 1}$ which is responsible for the fine probabilistic behavior of a lacunary sum. Takahashi~\cite{taka_1961} proved a CLT (with pure Gaussian limit) under the assumption that $n_{k+1}/n_k \to \infty$, and Gaposhkin~\cite{gap1966} proved that a CLT (with pure Gaussian limit) holds when $n_{k+1}/n_k$ is an integer for all $k$, or if $n_{k+1}/n_k \to \alpha$ for some $\alpha$ such that $\alpha^r \not\in \mathbb{Q},~r = 1,2,\dots$ (and if additionally the variance does not degenerate). A general framework connecting Diophantine equations and the distribution of lacunary sums was established in Gaposhkin's profound paper~\cite{gapo_1970}, where he proved that a CLT (with pure Gaussian limit) holds if for all fixed positive integers $a,b$ the number of solutions $k,\ell$ of the Diophantine equation \eqref{dio_equ} is bounded by a constant which is independent of $c$ (where only $c \neq 0$ needs to be considered, provided that the variance does not degenerate). One can check the validity of this general condition for sequences satisfying the assumptions mentioned earlier in this paragraph, such as $n_{k+1}/n_k \to \infty$ or $n_{k+1}/n_k \to \alpha$ for $\alpha^r \not\in \mathbb{Q}$. Finally, an optimal result was established in~\cite{abclt}: For $(n_k)_{k \geq 1}$ satisfying the Hadamard gap condition, the limit distribution of $N^{-1/2} \sum_{k=1}^N f(n_k x)$ is Gaussian provided that the number of solutions $(k,\ell)$ of \eqref{dio_equ}, subject to $k,\ell \leq N$, is of order $o(N)$ (for all fixed $a,b$, uniformly in $c \neq 0$). If, on the other hand, for some $a,b,c$ the number of solutions is $\Omega(N)$, then the CLT generally fails to hold. If the number of solutions with $c=0$ also is of order $o(N)$, then the CLT has the ``correct'' variance $\int_0^1 f(x)^2dx$, in perfect accordance with the independent case. Even if the number of solutions is of order $\Omega(N)$ for some $a,b,c$, then the deviation of the distribution of $N^{-1/2} \sum_{k=1}^N f(n_k x)$ from the Gaussian distribution can be quantified in terms of the ratio ``(number of solutions)/$N$''. This shows for example that while the CLT generally fails in the case $n_{k+1}/n_k \to p/q$, one obtains an ``almost CLT'' if both $p$ and $q$ are assumed to be large. Another example for such an ``almost CLT'' is when the growth constant in the Hadamard gap condition is assumed to be very large. See the statement of~\cite[Theorem 1.3]{abclt} and the subsequent discussion for more details.\\

Note that Gaposhkin's condition implies that the CLT also holds for all subsequences that are picked out of $(n_k)_{k \geq 1}$. This is not the case under the assumptions from~\cite{abclt}, where one might be able to extract a subsequence along which the CLT fails (by choosing a subsequence which allows a large number of solutions of the relevant Diophantine equations). It is interesting that the probabilistic behavior of lacunary sums might change when one passes to a subsequence of the original sequence -- this is in clear contrast to the bevahior of sums of independent random variables, where any subsequence of course is independent as well. A similar remark holds for permutations of lacunary sums resp.\ permutations of sums of independent random variables. These phenomena have received strong attention during the last years; see for example~\cite{abt_pams,abt_dep,abt_perm,a_e,fuk_perm}. To give only one sample result, in~\cite{abt_pams} the following is shown. As noted above, the CLT is true for pure trigonometric sums under the Erd\H os gap condition $n_{k+1}/n_k \geq 1 + c k^{-\alpha}$ for some $\alpha<1/2$. However, this is only true for the unpermuted sequence (i.e.\ sorted in increasing order). If permutations of the sequence are allowed, then this gap condition is not sufficient anymore for the validity of the CLT, as is no other gap conditon weaker than Hadamard's. More precisely, for any sequence $(\varepsilon_k)_{k \geq 1}$ with $\varepsilon_k \to 0$ there exists a sequence of positive integers satisfying $n_{k+1}/n_k \geq 1 + \varepsilon_k$, together with a permutation $\sigma: \mathbb{N} \mapsto \mathbb{N}$, such that the permuted (pure trigonometric) sum $N^{-1/2} \sum_{k=1}^N \cos (2 \pi n_{\sigma(k)} x)$ converges in distribution to a non-Gaussian limit. One can also construct such examples where the norming sequence $N^{-1/2}$ has to be replaced by $(\log N)^{1/2} N^{-1/2}$ and the limit is a Cauchy distribution, and examples where no limit distribution exists at all. See~\cite{abt_pams} for details on this particular result, and Chapter 3 of~\cite{be2017B} for a detailed discussion of permutation-invariance of limit theorems for lacunary (trigonometric) systems.\\

We close this section with some further references. For Hadamard lacunary $(n_k)_{k \geq 1}$, the limit distribution of $N^{1/2} D_N(n_k x)$ was calculated in~\cite{ab_lim}; under suitable Diophantine assumptions it coincides with the Kolmogorov distribution, which is  the distribution of the range of a Brownian bridge. A central limit theorem for Hardy--Littlewood--P\'olya sequences was established in~\cite{fuk_pet}. In~\cite{clb} the Erd\H os--Fortet example was revisited from the perspective of ergodic theory, and was interpreted in terms of the limiting behavior of certain modified ergodic sums, and generalized to cases such as expanding maps, group actions, and chaotic dynamical systems under the assumption of multiple decorrelation. See also~\cite{cilb,clb2}. The limit distribution of $N^{-1/2} \sum_{k=1}^N \cos (2\pi n_k x)$ for the special sequence $n_k=k^2,~k \geq 1,$ was determined by Jurkat and Van Horne in~\cite{jvh1,jvh2,jvh3}, and turned out to have finite moments of order $<4$, but not of order $4$. The theory of such sums is closely related to theta sums, and goes back to Hardy and Littlewood~\cite{hl}. For further related results, see~\cite{cema,fjk,walfisz}. For non-Gaussian limit distributions of $N^{-1/2} \sum_{k=1}^N \cos (2\pi n_k x)$ near the Erd\H{o}s gap condition $n_{k+1}/n_k \ge 1+ ck^{-1/2}$ see~\cite{berk_nong}. For a multidimensional generalization of Kac's results see~\cite{fuk_rr,fuk_rr2}, and for a multidimensional generalization of the CLT for Hardy--Littlewood--P\'olya sequences (considering a semi-group generated by powers of matrices instead) see~\cite{lev_2,lev}. See also~\cite{cc,cdv} for generalizations of the CLT for Hardy--Littlewood--P\'olya sequences to a very general setup of sums over powers of transformations/automorphisms.

\section{The law of the iterated logarithm for lacunary sequences} \label{sec:LIL}

Together with the law of large numbers (LLN) and the central limit theorem (CLT), the law of the iterated logarithm (LIL) is one of the fundamental results of probability theory. Very roughly speaking, the (strong) law of large numbers says that when scaling by $N^{-1}$ one has almost sure convergence of a sum of random variables, and the central limit theorem says that when scaling by $N^{-1/2}$ one has a (Gaussian) limit distribution. The law of the iterated logarithm operates between these two other asymptotic limit theorems; in its simplest form, it says that for a sequence $(X_n)_{n \geq 1}$ of centered i.i.d.\ random variables (under suitable extra assumptions, such as boundedness) one has
$$
\limsup_{N \to \infty} \frac{\sum_{n=1}^N X_n}{\sqrt{2 N \log \log N}} = \sigma \qquad \text{almost surely},
$$
where $\sigma$ is the standard deviation. Heuristically, the law of the iterated logarithm identifies the threshold between convergence in distribution and almost sure convergence for sums of i.i.d.\ random variables; indeed, while $\frac{\sum_{n=1}^N X_n}{\sqrt{2 N \log \log N}}$ converges to 0 in distribution by the CLT, it does not converge to 0 almost surely by the LIL. The first version of the LIL was given by Khinchin in 1924, and a more general variant was established by Kolmogorov in 1929. Note that the law of large numbers for trigonometric sums or sums of dilated functions is rather unproblematic: for \emph{any} sequence of distinct integers $(n_k)_{k \geq 1}$ one has
\begin{equation} \label{LIL}
\lim_{N \to \infty} \frac{1}{N} \sum_{k=1}^N f(n_k x) = \int_0^1 f(x)~dx,
\end{equation}
as long as one can assume a bit of regularity for $f$ (such as $f$ being a trigonometric polynomial, being Lipschitz-continuous, being of bounded variation on $[0,1]$, etc.). Only if one is not willing to impose any regularity assumptions upon $f$ the situation becomes quite different; see the remarks at the end of Section~\ref{sec:arithm}. \\

The most basic law of the iterated logarithm for lacunary systems is
\begin{equation} \label{LIL_nk}
\limsup_{N \to \infty} \frac{\sum_{n=1}^N \cos(2 \pi n_k x)}{\sqrt{2 N \log \log N}} = \frac{1}{\sqrt{2}}\qquad \text{a.e.}
\end{equation}
under the Hadamard gap condition on $(n_k)_{k \geq 1}$; this was obtained by Salem and Zygmund (upper bound)~\cite{salz} and Erd\H os and G\'al (lower bound)~\cite{egot}. Generally speaking, as often in probability theory the lower bound is more difficult to establish than the upper bound, since the latter can be proved by an application of the first Borel--Cantelli lemma (convergence part), while the former is proved by the second Borel--Cantelli lemma (divergence part, which needs some sort of stochastic independence as an extra assumption). Note that \eqref{LIL_nk} is a perfect analogue of \eqref{LIL} with the ``correct'' constant on the right-hand side.\\

As in the case of the CLT, replacing pure trigonometric sums by sums of more general 1-periodic functions makes the situation much more delicate. As in the previous section, a key role is played by Diophantine equations. However, while for the CLT it is crucial that the number of solutions of Diophantine equations ``stabilizes'' in some way to allow for a limit distribution (albeit a potentially non-Gaussian one), no such property is necessary for the validity of a form of the LIL (since, as noted above, this is defined as a $\limsup$, not as a $\lim$). Instructive examples are the following. In all examples, we assume that $f$ is 1-periodic with mean zero and bounded variation on $[0,1]$.
\begin{itemize}
 \item If $n_{k+1}/n_k \to \infty$ as $k \to \infty$, then
 $$
\limsup_{N \to \infty} \frac{\sum_{n=1}^N f(n_k x)}{\sqrt{2 N \log \log N}} = \left( \int_0^1 f^2(x)~dx\right)^{1/2} \qquad \text{a.e.}
 $$
\item If $n_{k} = 2^k,~ k \geq 1$, then
 $$
\limsup_{N \to \infty} \frac{\sum_{n=1}^N f(n_k x)}{\sqrt{2 N \log \log N}} = \sigma_f \qquad \text{a.e.},
$$
with
$$
\sigma_f^2 = \int_0^1 f^2(x)dx + 2 \sum_{m=1}^\infty \int_0^1 f(x) f(2^m x)~dx.
$$
\item Assume that $n_{k+1}/ n_k \geq q > 1,~k \geq 1$. Then there exists a constant $C$ (depending on $f$ and on $q$) such that
\begin{equation} \label{Tak_bound}
\limsup_{N \to \infty} \frac{\sum_{n=1}^N f(n_k x)}{\sqrt{2 N \log \log N}} \leq C \qquad \text{a.e.}
\end{equation}
\item If $n_k = 2^k-1,~k \geq 1,$ and if $f(x) = \cos(2 \pi x) + \cos(4 \pi x)$, then
\begin{equation} \label{kac_non}
\limsup_{N \to \infty} \frac{\sum_{n=1}^N f(n_k x)}{\sqrt{2 N \log \log N}} = \sqrt{2} |\cos(\pi x)| \qquad \text{a.e.}
\end{equation}
\end{itemize}

The first result in this list (due to Takahashi~\cite{tak_inf}) is in perfect accordance with the LIL for truly independent random sums, in accordance with the fact that also the CLT holds in the ``truly independent'' form under the large gap condition $n_{k+1}/n_k \to \infty$. The second result is an analogue of Kac's CLT in Equations \eqref{kac1a} and \eqref{kac2a}: as with the CLT, also the LIL holds for the sequence $(2^k)_{k \geq 1}$, but the limiting variance deviates from the one in the ``truly independent'' case. Note that in contrast to the CLT case we now do not need to require that $\sigma_f \neq 0$ for the validity of the statement. The third result (Takahashi~\cite{tak_had}) asserts that there is an upper-bound version of the LIL for lacunary sums (even for sequences where there is no convergence of distributions, and any form of the CLT fails). Finally, the fourth result (the Erd\H os--Fortet example for the LIL instead of the CLT) shows the remarkable fact that the $\limsup$ in the LIL for Hadamard lacunary sums might actually be \emph{non-constant} -- this is very remarkable, and a drastic deviation from what one can typically observe for sequences of independent random variables. In particular this example shows that under the Hadamard gap condition an upper-bound version of the LIL is in general the best that one can hope for. Not very surprisingly, the source of all these phenomena are (as in the previous section) Diophantine equations such as \eqref{dio_equ}, and their number of solutions within the sequence $(n_k)_{k \geq 1}$. So in the LIL there is again a complex interplay between probabilistic, analytic and arithmetic aspects which controls the fine asymptotic behavior of lacunary sums.\\

In probability theory there is a version of the LIL for the Kolmogorov--Smirnov statistic of an empirical distribution. This is called the Chung--Smirnov LIL, and in the special case of a sequence $(X_n)_{n \geq 1}$ of i.i.d.\ random variables having uniform distribution on $[0,1]$ (where the Kolmogorov--Smirnov statistic coincides with the discrepancy) it asserts that
$$
\limsup_{N \to \infty} \frac{N D_N(X_n)}{\sqrt{2 N \log \log N}} = \frac{1}{2} \qquad \text{almost surely.}
$$
Here the number $1/2$ on the right-hand side arises essentially as the maximal $L_2$ norm (``standard deviation'') of a centered indicator function of an interval $A \subset [0,1]$ (namely the indicator function of an interval of length $1/2$). Based on the principle that lacunary sequences tend to ``imitate'' the behavior of truly independent sequences, it was conjectured that an analogue of the Chung--Smirnov LIL should also hold for the discrepancy of $(\{n_k x\})_{k \geq 1}$, where $(n_k)_{k \geq 1}$ is a Hadamard lacunary sequence. This was known as the Erd\H os--G\'al conjecture, and was finally solved by Philipp~\cite{plt}, who proved that for any $q>1$ there exists a constant $C_q$ such that for $(n_k)_{k \geq 1}$ satisfying $n_{k+1}/n_k \geq q$ we have
\begin{equation} \label{Phil_LIL}
\frac{1}{\sqrt{32}} \leq \limsup_{N \to \infty} \frac{N D_N(\{n_k x\})}{\sqrt{2 N \log \log N}} \leq C_q \qquad \text{a.e}.
\end{equation}
An admissible value of $C_q$ was specified in~\cite{plt} as $C_q = 166/\sqrt{2} + 664/(\sqrt{2q} - \sqrt{2})$. The first inequality in \eqref{Phil_LIL} follows from (a complex version of) Koksma's inequality together with \eqref{LIL_nk}, so the novelty is the second inequality (upper bound). Note also that the upper bound in \eqref{Phil_LIL} implies Takahashi's ``upper bound'' LIL in \eqref{Tak_bound}, again as a consequence of Koksma's inequality.\\

Philipp's result has been extended and refined into many different directions. The most precise results, many of which were obtain by Fukuyama, show again a fascinating interplay between arithmetic, analytic and probabilistic effects. As a sample we state the following results (all from Fukuyama's paper~\cite{fuk_theta}):\\

Let $n_k = \theta^k,~k  \geq 1$. Then
\begin{equation} \label{LIL_dn}
\limsup_{N \to \infty} \frac{N D_N(\{n_k x\})}{\sqrt{2 N \log \log N}}
 \end{equation}
 exists and is constant (almost everywhere). Denoting the value of this $\limsup$ by $\Sigma_\theta$, we have:
\begin{itemize}
\item If $\theta^r \not\in \mathbb{Q}$ for all $r = 1,2,\dots$, then $\Sigma_\theta = 1/2$ a.e.
\item If $r$ denotes the smallest positive integer such that $\theta^r=p/q$ for some coprime $p,q$, then $1/2 \leq \Sigma_\theta \leq \sqrt{(pq+1)/(pq-1)}/2$ a.e.
\item If $\theta^r = p/q$ as above and both $p$ and $q$ are odd, then $\Sigma_\theta = \sqrt{(pq+1)/(pq-1)}/2$ a.e.
\item If $\theta=2$, then $\Sigma_\theta = \sqrt{42}/9$ a.e.
\item If $\theta > 2$ is an even integer, then $\Sigma_\theta = \sqrt{(p+1)p(p-2)/(p-1)^3}/2$ a.e.
\item If $\theta = 5/2$, then $\Sigma_\theta = \sqrt{22}/9$ a.e.
\end{itemize}
All these results were obtained by very delicate calculations involving Fourier analysis and Diophantine equations. The calculations from~\cite{fuk_theta} were continued by the same author and his group in~\cite{fuk_c3,fuk_c1,fuk_c2,fuk_c4}, so that now we have a relatively comprehensive picture on the behavior of these $\limsup$'s in the case when $(n_k)_{k \geq 1}$ is (exactly) a geometric progression.\\

In~\cite{a_lil,a_lil2} for general Hadamard lacunary sequences $(n_k)_{k \geq 1}$ a direct connection was established which links the number of solutions of \eqref{dio_equ} with the value of the $\limsup$ in the LIL, in the same spirit as this was done before in~\cite{abclt} for the CLT (as described in the previous section). In particular, if the number of solutions of \eqref{dio_equ} is sufficiently small, then the LIL holds with the constant $1/2$ on the right-hand side, exactly as in the truly independent case. Another interesting observation is that if $(n_k)_{k \geq 1}$ is Hadamard lacunary with growth factor $q>1$, and if $\Sigma$ denotes the value of the $\limsup$ in \eqref{LIL_dn}, then the difference $|\Sigma-1/2|$ can be quantified in terms of $q$ and tends to zero a.e.\ as $q \to \infty$. Thus there is a smooth transition towards the ``truly independent'' LIL as the growth factor $q$ increases, and under the large gap condition $n_{k+1}/n_k \to \infty$ the value of $\Sigma$ actually equals 1/2. Another remarkable fact is that there exist Hadamard lacunary sequences for which the $\limsup$ in the LIL for the discrepancy is not a constant almost everywhere, but rather a function of $x$, similar to what happened in \eqref{kac_non} for the LIL for $\sum f(n_k x)$. In some cases the limit functions in the LIL for the discrepancy can be explicitly calculated, and are ``surprisingly exotic'' (in the words of Ben Green's MathSciNet review of~\cite{aist_irr2}).\\

As noted above, Philipp's LIL for the discrepancy has been extended into many different directions. For example, while it is known that the result can fail as soon as the Hadamard gap condition is relaxed to any sub-exponential growth condition, it turns out to be possible to obtain an LIL for the discrepancy when a weaker growth condition is compensated by stronger arithmetic assumptions. In particular, an analogue of Philipp's result has been proved for Hardy--Littlewood--P\'olya sequences~\cite{phl_emp}; see also~\cite{a_sub,bpt_emp,fuk_nak,tech_zaf}. As a closing remark concerning the LIL, it is interesting that the optimal value of the lower bound in \eqref{Phil_LIL} is still unknown; cf.\ \cite{aff} for more context.\\

While much effort has been spent towards understanding the probabilistic behavior of lacunary sums at the scales of the CLT and LIL, it seems that investigations at other scales (such as in particular at the large deviations scale) were only started recently. The few results which are currently available point once again towards an intricate connection between probabilistic, analytic and arithmetic effects; see the very recent papers \cite{AGKPR,FJP}.

\section{Normality and pseudorandomness} \label{sec:norm}

Normal numbers were introduced by Borel~\cite{borel} in 1909.  From the very beginning the concept of normality of real numbers was associated with ``randomness''. While normality of real numbers was originally defined in terms of counting the number of blocks of digits, it is not difficult to see\footnote{Probably first explicitly mentioned by D.D.\ Wall in his PhD thesis, 1949.} that a number $x$ is normal in base $b$ if and only if the sequence $(\{ b^n x \})_{n \geq 1}$ is equidistributed. As Borel proved, Lebesgue-almost all real numbers are normal in a (fixed) integer base $b \geq 2$, and thus almost all reals are normal in \emph{all} bases $b \geq 2$ (such numbers are called \emph{absolutely normal}). While normal numbers are ubiquitous from a measure-theoretic perspective,\footnote{Interestingly, while the set of normal numbers is large from a measure-theoretic point of view, it turns out to be small from a topological point of view. More precisely, the set of normal numbers is meager (of first Baire category), see e.g.\ \cite{gsw}. In the words of Edmund Hlawka~\cite[p.~78]{hlawka}: ``Thus whereas the normal numbers almost force themselves on to the measure theorist, the topologist is apt to overlook them entirely.''} it is difficult to construct normal numbers. The most fundamental construction is due to Champernowne, who proved (using combinatorial arguments) that the number
$$
0. 1~2~3~4~5~6~7~8~9~10~11~12~13~14~\dots,
$$
which is obtained by a concatenation of the decimal expansions of the integers, is normal in base 10. The idea of creating normal numbers by a concatenation of the ($b$-ary) expansions of the values of (simple) functions at integers (or primes) is still the most popular, and probably most powerful, method in this field. We only note that Copeland and Erd\H os~\cite{coperd} proved that
$$
0.2~3~5~7~11~13~17~19~23~\dots,
$$
which is obtained by concatenating the decimal expansions of the primes, is normal in base 10, and refer to~\cite{daver,dk,mad,mtt,nak1,nak2} for more results of this flavor. It should be noted, however, that there have been earlier constructions of a conceptually very different nature, such as that of Sierpinski~\cite{sier} in 1917. See~\cite{bfsier} for an exposition of Sierpinski's construction and more context, and see also~\cite{bfs_tur} on an early (unpublished) algorithm of Turing for the construction of normal numbers. We finally mention a very recent idea for the construction of a normal number by Drmota, Mauduit and Rivat \cite{DMR}, which is not based on the concatenation of decimal blocks as above, but rather on the evaluation of an automatic sequence along a subsequence of the index set (in this particular case, the Thue--Morse sequences evaluated along the squares); see also \cite{Mspiegel,spiegel}.\\

While Lebesgue-almost all numbers are normal and there are some constructions of normal numbers, it is generally considered to be completely hopeless to prove that natural constants such as $\pi,e,\sqrt{2},\dots$ are normal in a given base (although the experimental evidence clearly points in that direction:~\cite{nspi,wagon,xyl}). Many such open problems ``for the next millennium'' are contained in Harman's survey article~\cite{harman100}; see also~\cite{pi_day}. However, it is quite clear that the mathematical machinery which would be necessary to prove the normality of $\sqrt{2}$ or other such constants is completely lacking; compare the rather deplorable current state of knowledge on the binary digits of $\sqrt{2}$ as given in~\cite{BBCP,dub,van}. A small spark of hope is provided by the very remarkable formulas of Bailey, Borwein and Plouffe (now widely known as BBP formulas), which allow to calculate deep digits of $\pi$ (and other constants) without the need of computing all previous digits. See~\cite{pi_source} for a very comprehensive ``source book'' covering computational aspects of $\pi$, and see~\cite{bailey_c} for a very rare example of a possible strategy of what a proof of the normality of $\pi$ could possibly like (cf.\ also~\cite{lag}).\\

Since normality of $x$ in a base $b$ can be expressed in terms of the equidistribution of the sequence $(\{b^n x\})_{n \geq 1}$, it is very natural to consider the discrepancy of $D_N (\{b^n x\})$ and call this (with a slight abuse of language) the discrepancy of $x$ (as a normal number, with respect to a base $b$). Remarkably, it is still unknown how small the discrepancy of a normal number can be (this is known as \emph{Korobov's problem}). Levin~\cite{levin_n} constructed (for given base $b$) a number $x$ such that
$$
D_N(\{b^n x\}) = O \left(\frac{(\log N)^2}{N} \right);
$$
by Schmidt's general lower bound the exponent of the logarithm cannot be reduced below 1, but the optimal size of this exponent remains open.\\

One of the most interesting, and most difficult, aspects of normal numbers is normality with respect to two or more different bases. Extending work of Cassels~\cite{cass_norm}, Schmidt~\cite{schm} characterized when normality with respect to a certain base implies normality with respect to another base, and when this is not the case. See also~\cite{bso,bug}. However, generally speaking it is very difficult to construct numbers which are normal with respect to several different bases, and the ``constructions'' are much less explicit than the ones of Champernowne and Copeland--Erd\H os mentioned above.  The problem of the minimal order of the discrepancy of normal numbers seems to be very difficult when different bases are considered simultaneously. Aistleitner, Becher, Scheerer and Slaman~\cite{abss} constructed a number $x$ such that
$$
D_N(\{b^n x\}) = O_b \left(N^{-1/2} \right)
$$
for all integer bases $b \geq 2$; this is considered to be an ``unexpectedly small'' order of the discrepancy by Bugeaud in his MathSciNet review of~\cite{abss}. It is not known if the exponent $-1/2$ of $N$ in this estimate is optimal or not; indeed, no non-trivial lower bounds whatsoever (beyond the general lower bound $(\log N) / N$ of Schmidt) are known for this problem, but it is quite possible that whenever simultaneous normality with respect to different (multiplicatively independent) bases is considered, there must be at least one base for which the discrepancy is ``large''.\\

In a recent years, there has been a special focus on algorithmic aspects of the construction of normal numbers. A particularly striking contribution was a polynomial-time algorithm for the construction of absolutely normal numbers due to Becher, Heiber and Slaman~\cite{bhs}. See also~\cite{alfig, beli, schee}. Related to such algorithmic and computational problems are questions on the complexity of the set of normal numbers from the viewpoint of descriptive set theory in mathematical logic; in this framework, the rank of the set of normal numbers~\cite{kili} and absolutely normal numbers~\cite{behe} within the Borel hierarchy has been determined.\\

The notion of normality can be extended in a natural way to many other situations, where it is always understood that normality should be the typical behaviour. For example, one can consider normal continued fractions, where the ``expected'' number of occurences of each partial quotient is prescribed by the Gauss--Kuzmin measure; see for example~\cite{akm,by}. Other generalizations consider for example normality with respect to $\beta$-expansions~\cite{bk,mst}, a numeration system which generalizes the $b$-ary expansion to non-integral bases $\beta$, or normality with respect to Cantor expansions~\cite{amv,fs,man}, a numeration system which allows a different set of ``digits'' at each position. For a particularly general framework, see~\cite{mm}. Interestingly, in such generalized numeration systems there can be more than one natural definition of normality, using as starting point for exampe either the idea of counting blocks of digits, or the idea of equidistribution of an associated system. The relation between such different (sometimes conflicting) notions of normality has been studied in particular detail for Cantor expansions~\cite{airm,mance2,mance1}.   \\

Normal numbers feature prominently in the chapter on random numbers in Volume~2 of Knuth's celebrated series on \emph{The Art of Computer Programming}~\cite{knuth}. There he tries to come to terms with the notion of ``random'' sequences of numbers, and introduces an increasingly restrictive scheme of ``randomness'' of deterministic sequences. The concept of normality is also the starting point for one of the (quantitative) measures of pseudorandomness, which were introduced by Mauduit and S\'{a}rk\"{o}zy~\cite{maus} and then studied in a series of papers. Note in this context that the transformation $T:~x \mapsto b x$ mod 1, which is at the foundation of the concept of normal numbers, can in some sense be seen as the continuous analogue of the recursive formula which defines a linear congruential generator (LCG), one of the most classical devices for pseudorandom number generation. For this connection between normal numbers and pseudo-random number generators, see for example~\cite{BCrand}. Another very fruitful aspect of normal numbers is the connection with ergodic theory, which comes from the observation that the sequence $(\{b^n x \})_{n \geq 1}$ is the orbit of $x$ under the transformation $T$ from above, and that this transformation is measure-preserving (with respect to the Lebesgue measure) and ergodic. We will not touch upon this connection in any detail, and instead refer to~\cite{da_ka,ke_sa}.\\

Another sequence which is often associated with ``randomness'' is the sequence $(\{x^n\})_{n \geq 1}$ for real $x>1$, or more generally $(\{\xi x^n\})_{n \geq 1}$ for $\xi \neq 0$ and $x>1$. This looks quite similar to a (Hadamard) lacunary sequence such as $(b^n x)_{n \geq 1}$, but its metric theory is of a very different nature in several respects. Both sequences are variants of a geometric progression, but while in the lacunary sequence the base $b$ is fixed and $x$ is assumed to be a ``parameter'', now $\xi$ is assumed to be fixed and the base $x$ of the geometric progression is the parameter. While $(b^n x)_{n \geq 1}$ can in many ways be easily interpreted in terms of harmonic analysis, digital expansions, ergodic theory, etc., such simple interpretations fail for $(\{x^n\})_{n \geq 1}$. Note in particular that in contrast to lacunary sequences there now is no periodicity when replacing $x \mapsto x+1$, there is no ``orthogonality'', and the calculation of moments of sums $\sum f(x^n)$ does not simply reduce to the counting of solutions of Diophantine equations. Still, what is preserved from the setup of lacunary sequences is that $x^n$ (as a function of $x$) oscillated quickly on intervals where $x^m$ is essentially constant, provided that $n$ is significantly larger than $m$, and there are good reasons to consider systems such as $(\cos(2 \pi x^n))_{n \geq 1}$ to be ``quasi-orthogonal'' and ``almost independent'' in some appropriate sense.\\

One of the most fundamental results on this type of sequence is due to Koksma~\cite{koks_x}: assuming that $\xi \neq 0$ is fixed, the sequence $(\{\xi x^n\})_{n \geq 1}$ is uniformly distributed mod 1 for almost all $x>1$. In particular, when $\xi=1$, the sequence $(\{x^n\})_{n \geq 1}$ is uniformly distributed mod 1 for almost all $x > 1$. In very sharp contrast with Koksma's metric result is the fact that until today \emph{not a single example} of a number $x$ is known for which $(\{x^n\})_{n \geq 1}$ is indeed uniformly distributed. This problem is related with Mahler's problem on the range of $(\{(3/2)^n\})_{n \geq 1}$, which also seems to be completely out of reach for current methods (cf.\ \cite{dubi,flp}). The sequence $(\{x^n\})_{n \geq 1}$ is discussed at length in Knuth's book, where it is conjectured that this sequence is a good candidate to pass several very strict pseudorandomness criteria for almost all $x$. For example, Knuth conjectured that for all sequences of distinct integers $(s_n)_{n \geq 1}$ the sequence $(\{x^{s_n}\})_{n \geq 1}$ (a subsequence of the original sequence) has a strong equidistribution property called complete uniform distribution, for almost all $x > 1$; this was indeed established by Niederreiter and Tichy~\cite{niti}. It is also known that $(\{x^n\})_{n \geq 1}$ satisfies an law of the iterated logarithm in the ``truly independent'' form
$$
\limsup_{N \to \infty} \frac{N D_N(\{x^n\})}{\sqrt{2 N \log \log N}} = \frac{1}{2} \qquad \text{a.e}.,
$$
and similarly satisfies a central limit theorem which is perfectly analogous to the one for truly independent systems~\cite{a_g}. Note that Knuth's assertion that the sequence $(\{x^n\})_{n \geq 1}$ shows good pseudo-random behavior for almost all $x>1$ is of limited practical use, as long as no such value of $x$ is found. The discrete analogue would be to study the pseudo-randomness properties of $a^n$ mod $q$ for $n=1,2,\dots$, where $a$ and $q$ are fixed integers. Investigations on the pseudo-randomness properties of such sequence were carried out for example by Arnol'd~\cite{arnold}, who experimentally observed good pseudorandom behavior; cf. also~\cite{aist_ar}.\\

To close this section, we note that equidistribution is of course just one property which can be used to characterize ``pseudorandom'' behavior (essentially by analogy with the Glivenko--Cantelli theorem). There are many other statistics which could be applied to a sequence in $[0,1]$ to determine whether it behaves in a ``random'' way or not. One class of such statistics are gap statistics at the level of the average gap (which is of order $1/N$ when considering the first $N$ elements of a sequence in $[0,1]$), such as the distribution of nearest-neighbor gaps, or the pair correlation statistics. We do not give formal definitions of these concepts here, but note that they are inspired by investigations of the statistics of quantum energy eigenvalues in the context of the Berry--Tabor conjecture in theoretical physics; see~\cite{mark} for more context. Pseudorandom behavior with respect to such statistics is called ``Poissonian'', since it agress with the corresponding statistics for the Poisson process. The general principle that lacunary systems show pseudorandom behavior is also valid in this context. For example, Rudnick and Zaharescu~\cite{rud_z} showed that for $(n_k)_{k \geq 1}$ satisfying the Hadamard gap condition the sequence $(\{n_k x\})_{k \geq 1}$ is Poissonian for almost all $x$, and Aistleitner, Baker, Technau and Yesha~\cite{aist_arx} showed that the same holds for $(\{x^n\})_{n \geq 1}$ for almost all $x>1$.\\

This section on normal numbers and sequences of the form $(\{\xi x^n\})_{n \geq 1}$ gives of course only a very brief overview of the subject, and has to leave out many interesting aspects. For a much more detailed exposition we refer the reader to the book of Bugeaud \cite{bug_book}.

\section{Random sequences} \label{sec:rand}

In the previous sections we have illustrated the philosophy that gap sequences exhibit many probabilistic properties which are typical for sequences of i.i.d.\ random variables. In many cases the large gap condition $n_{k+1}/n_k \to \infty$ gives ``true'' random limit theorems, the Hadamard gap condition $n_{k+1}/n_k \geq q > 1$ is a critical transition point where a mixture of probabilistic, analytic and arithmetic effects comes into play, and the ``almost independent'' behavior is lost when the gap condition is relaxed below Hadamard's. There are results which hold under weaker gap conditions such as the Erd\H os gap condition $n_{k+1}/n_k \geq 1 + c k^{-\alpha}$, $0<\alpha<1/2$, or under additional arithmetic assumptions, but as a whole the Hadamard gap condition is the critical point where the ``almost independent'' behavior of systems of dilated functions starts to break down.\\

However, while almost independent behavior is generally lost under a weaker gap condition (without strong arithmetic information), there is another possible perspective on the problem. As noted, for a fixed sequence $(n_k)_{k \geq 1}$ one cannot expect ``almost independent'' behavior of $(\{n_k x\})_{k \geq 1}$, say, without assuming a strong growth condition on $(n_k)_{k \geq 1}$. However, even without such a growth condition one might expect that $(\{n_k x\})_{k \geq 1}$ shows independent behavior for ``typical'' sequences $(n_k)_{k \geq 1}$. Here the word ``typical'' of course implies that the sequence has to be taken from a generic set in some appropriate space which possesses a measure, so quite naturally this idea leads to considering ``random'' sequences $(n_k)_{k \geq 1} = (n_k(\omega))_{k \geq 1}$ which are constructed in a randomized way over some probability space.\\

Of course there are many possible ways how a random sequence can be constructed. From results of Salem and Zygmund~\cite{sal_z} for trigonometric sums
with random signs it follows easily that if we define a sequence $(n_k)_{k \geq 1}$ by flipping a coin (independently) for every positive integer to decide whether it should be contained in the sequence or not, and let $\mathbb{P}$ denote the probability measure on the space over which the ``coins'' are defined, then for $\mathbb{P}$-almost all sequences as defined above one has
\begin{equation}\label{CLTnew}
\frac{1}{\sqrt N} \sum_{k=1}^N \cos (2\pi n_kx) \overset{\mathcal D} {\longrightarrow}
N(0,1/4)
\end{equation}
and
\begin{equation}\label{LILnew}
\limsup_{N\to\infty} \, \frac{1}{\sqrt{2N \log\log N}} \sum_{k=1}^N
\cos (2\pi n_kx) = \frac{1}{2} \qquad \text{for almost all $x$},
\end{equation}
where $N(0,\sigma^2)$ denotes the normal distribution with mean 0 and variance $\sigma^2$ and $\overset{\mathcal D} {\longrightarrow}$ denotes convergence in distribution. Note that (\ref{CLTnew}) and (\ref{LILnew}) are  not exactly matching with the truly independent case, where the limit distribution would be $N(0, 1/2)$ and the limsup in the LIL would be $1/\sqrt{2}$. The ``loss'' on the right-hand sides of \eqref{CLTnew} and \eqref{LILnew} comes from the fact that
a Dirichlet kernel is ``hiding'' in this linearly growing sequence, and this kernel is highly localized near 0 and 1 so that its contribution is lost in the CLT and LIL. By the strong law of large numbers (SLLN) clearly $n_k \sim 2k$ as $k \to \infty$, $\mathbb{P}$-almost surely, so the sequences constructed here are very far from satisfying any substantial gap condition; in contrast, their (typical) order of growth is only linear. It should be noted that the gaps $n_{k+1}-n_k$ in this sequence are not bounded: with full $\mathbb P$-probability, $n_{k+1}-n_k=1$ for infinitely many $k$ (roughly, in half of the cases), but for infinitely many $k$, the gap $n_{k+1}-n_k$ has order of magnitude $c\log k$; this follows from the ``pure heads'' theorem of Erd\H{o}s and R\'enyi, see~\cite{renyi70}.\\

We call an increasing sequence $(n_k)_{k\ge 1}$ of positive integers a $B_2$ sequence if there exists a constant $C>0$ such that for any integer
$\nu>0$ the number of representations of $\nu$ in the form $\nu=n_k\pm n_\ell,~k>\ell \ge 1$, is at most $C$. By a result of Gaposhkin~\cite{gapo_1970} already mentioned in Section~\ref{sec:CLT}, the sequence $(f(n_k x))_{k\ge 1}$ satisfies the CLT for all Hadamard lacunary $(n_k)_{k\ge 1}$ and all $1$-periodic Lipschitz continuous $f$ if and only if for any $m\ge 1$, the set-theoretic union of the sequences $(n_k)_{\ge 1}$,  $(2n_k)_{\ge 1}$, \ldots, $(m n_k)_{\ge 1}$ satisfies the $B_2$ condition.\footnote{Note that the definition of the $B_2$ property used in~\cite{gapo_1970} is slightly different from the standard usage in  number theory (see e.g.\ \cite{halb}) requiring that the number of solutions of $\nu=n_k + n_\ell, ~k>\ell \ge 1$, is bounded by $C$, but this does not affect the discussion below.}
No similarly complete result is known for sequences $(n_k)_{k\ge 1}$ growing slower than exponentially, but Berkes~\cite{berk_rand0} proved that if $(n_k)_{k\ge 1}$ is a $B_2$ sequence satisfying the gap condition
\begin{equation} \label{erdgap}
n_{k+1}/n_k \ge 1+c k^{-\alpha}, \qquad k \geq 1,
\end{equation}
for some $c>0$, $\alpha>0$, then $(\cos (2\pi n_kx))_{k \geq 1}$ satisfies the CLT and LIL.
To verify the $B_2$ property for a concrete sequence $(n_k)_{k\ge 1}$ is generally a difficult problem, but the situation is quite different for random constructions. Let $I_1, I_2, \ldots$ be disjoint blocks of consecutive integers and let $n_1, n_2, \ldots$ be independent random variables on some probability space $(\Omega, \mathcal{F}, \mathbb{P})$ such that $n_k$ is uniformly distributed over $I_k$. Clearly, the number of different sums $\pm n_{k_1}\pm n_{k_2} \pm n_{k_3}$, $1\le k_1, k_2, k_3\le k-1,$ is at most $8(k-1)^3$, and thus if the size of $I_k$ is $\ge k^5$, then the probability that $n_k$ is equal to any of these sums is $\le 8k^3 k^{-5}=O(k^{-2})$. Thus by the Borel-Cantelli lemma, with $\mathbb P$-probability 1, such a coincidence can occur only for finitely many $k$. Thus the equation
$$\pm n_{k_1}\pm n_{k_2} \pm n_{k_3}\pm n_{k_4}=0,  \quad k_1\le k_2 \le k_3 <k_4$$
has only finitely many solutions, which implies that $(n_k)_{k \geq 1}$ is a $B_2$ sequence. \\

Let us recall now that by a result of Erd\H{o}s ~\cite{erd_o}, $(\cos (2\pi n_kx))_{k \geq 1}$ satisfies the CLT with limit distribution $N(0, 1/2)$, provided that (\ref{erdgap}) holds with $\alpha<1/2$, and this result is sharp, i.e.\ there exists a sequence $(n_k)_{k\ge 1}$ satisfying (\ref{erdgap}) with $\alpha=1/2$ such that the CLT fails for $(\cos (2\pi n_kx))_{k \geq 1}$. Note that the counterexample is irregular: while $n_{k+1}/n_k-1$ is of the order $O(k^{-1/2})$ for most $k$, there is also a subsequence along which $n_{k+1}/n_k\to\infty$.  One may therefore wonder if regular behavior of $n_{k+1}/n_k$ implies the CLT; in particular, Erd\H{o}s~\cite{erd_o} conjectured that the CLT holds for $(\cos (2\pi n_kx))_{k \geq 1}$ if $n_k= \lfloor e^{(k^\beta)} \rfloor$ for some $\beta$ in the range $0<\beta\le 1/2$. (Note that for $\beta>1/2$ condition (\ref{erdgap}) is satisfied with $\alpha<1/2$, so the CLT follows from Erd\H{o}s' result.)  This conjecture was proved by Murai~\cite{murai} for $\beta>4/9$, but for smaller $\beta$ the problem is still open. Random constructions provide here important information. Kaufman~\cite{kauf} proved that if $c$ is chosen at random, with uniform distribution on a finite interval $(a, b)\subset (0, \infty)$, then $(\cos (2\pi n_kx))_{k \geq 1}$ with $n_k=e^{(ck^\beta)}$ satisfies the CLT with probability 1 for any fixed $\beta>0$. An even wider class of random sequences with the CLT property is obtained by choosing the blocks $I_k$
in the random construction above as the integers in the interval
\begin{equation} \label{jk}
J_k= \left(e^{(ck^\beta)}(1-r_k), \,  e^{(ck^\beta)}(1+r_k)\right), \qquad r_k=o(k^{-(1-\beta)}).
\end{equation}
A simple calculation shows that these intervals are disjoint for $k\ge k_0$ and for $n_k\in J_k$ we have (\ref{erdgap}) with $\alpha=1-\beta$, in fact we even have
$$ n_{k+1}/n_k=1+ c_1(1+o(1))/k^{1-\beta} $$
with some constant $c_1>0$.
Now if $r_k$ decreases like a negative power of $k$, then the length of $J_k$ will be $\ge k^5$ and thus the constructed random sequence $(n_k)_{k \geq 1}$ will be a $B_2$ sequence with probability 1, so $(\cos (2\pi n_kx))_{k \geq 1}$ satisfies the CLT.  In other words, the CLT for $(\cos (2\pi n_kx))_{k \geq 1}$ holds for a huge class of sequences $n_k\sim e^{(ck^\beta)}$ for any $c>0$, $\beta>0$.\\

Concerning $B_2$ sequences, it is worth pointing out that Erd\H{o}s~\cite{er_b2} proved, decades before Carleson's convergence theorem, that $\sum_{k=1}^\infty (a_k \cos (2\pi n_kx) + b_k \sin (2\pi n_kx))$ is almost everywhere convergent if $(n_k)_{k\ge 1}$ is a $B_2$ sequence.
The question of how slowly $B_2$ sequences can grow is a much investigated problem of number theory, see e.g.\ Halberstam and Roth~\cite{halb}, Chapters II and III. It is easily seen that a $B_2$ sequence $(n_k)_{k\ge 1}$ cannot be $o(k^2)$ and Erd\H{o}s and R\'enyi~\cite{erb2} proved by a random construction that for any $\ve>0$ there exists a $B_2$ sequence $(n_k)_{k\ge 1}$ with $n_k = O(k^{2+\ve})$. Changing the $B_2$ property slightly and requiring that all numbers $n_k\pm n_\ell$, $k >\ell$, are actually different, makes the problem considerably harder. The  ``greedy algorithm'' yields a $B_2$ sequence $(n_k)_{k\ge 1}$ with $n_k=O(k^3)$, see~\cite{miac}, and it took nearly 40 years to improve this to $n_k=o(k^3)$, see~\cite{aksz}.  The best currently known (random) construction is due to Ruzsa~\cite{ruzsa}, and satisfies $n_k= k^{1/(\sqrt{2}-1)+o(1)}$.\\

Let $(\omega_n)_{n\ge 1}$ be a nondecreasing sequence of positive integers tending to $+\infty$ and let us divide the set of positive integers into disjoint blocks $I_1, I_2, \ldots$ such that the cardinality of $I_k$ is $\omega_k$. Using these blocks in the random construction above, the resulting random sequence $(n_k)_{k\ge 1}$ cannot be a $B_2$ sequence if $(\omega_n)_{n\ge 1}$ grows slower than any power of $n$, but it is proved in Berkes~\cite{berk_rand} that with $\mathbb P$-probability 1, $(\cos (2\pi n_k x))_{k \geq 1}$ still satisfies the CLT and LIL. The limit distribution here is $N(0, 1/2)$ and the limsup in the LIL is $1/\sqrt{2}$, so that the ``loss of mass'' phenomenon observed in the case of the random sequence $(n_k)_{k\ge 1}$ in the Salem-Zygmund paper~\cite{sal_z} does not occur here. The gaps in this sequence satisfy $n_{k+1}-n_k \le 2\omega_{k+1}$, i.e.\ they can grow arbitrarily slowly.  An LIL for the discrepancy of $(\{n_kx\})_{k \geq 1}$ under the same gap condition was given in Fukuyama~\cite{fuku_rand}. In~\cite{berk_rand} the question was raised if there exists a sequence $(n_k)_{k\ge 1}$ with bounded gaps $n_{k+1}-n_k=O(1)$ such that the CLT holds. Bobkov and G\"otze~\cite{bob_g}  showed that if we want no loss of mass in the CLT, the answer is negative: if $(n_k)_{k\ge 1}$ is any increasing sequence of positive integers with $n_{k+1}-n_k\le L,~k \geq 1,$ such that $N^{-1/2} \sum_{k=1}^N \cos (2\pi n_kx)$ has a Gaussian limit distribution $N(0, \sigma^2)$, then necessarily $\sigma^2 <1/2$ and $L\ge 1/(1-2\sigma^2)$. On the other hand, Fukuyama~\cite{fuk_ph,fuku_rand3,fuku_rand2} showed that for  any $\sigma^2<1/2$ there exists indeed a random subsequence $(\cos (2\pi n_kx))_{k \geq 1}$  of the trigonometric system satisfying the CLT with limit distribution $N(0, \sigma^2)$ and with bounded gaps $n_{k+1}-n_k\le L$ with $L \sim 4/(1-2\sigma^2)$ as $\sigma^2\to 1/2$.  This shows that the result of Bobkov and G\"otze is optimal up to a factor 4. This remarkable result is the ``small gaps'' counterpart of Erd\H{o}s' central limit theorem~\cite{erd_o}:  the latter determines the smallest gap sizes in $(n_k)_{k\ge 1}$  implying the CLT for $(\cos (2\pi n_kx))_{k \geq 1}$, while Fukuyama's result determines the smallest gap size which still allows a CLT with limit distribution
$N(0, \sigma^2)$ to hold.\\

It is worth pointing out that the bounded gap sequences  in~\cite{fuk_ph,fuku_rand3,fuku_rand2} are obtained by rather complicated random constructions, while using the previously discussed simple construction and choosing the $n_k$ as independent random variables uniformly distributed over adjoining  blocks $I_k$ with equal length results in a random sequence $(n_k)_{k\ge 1}$ satisfying almost surely
\begin{equation}\label{CLTBG}
\frac{1}{\sqrt N} \sum_{k=1}^N \cos (2\pi n_kx) \overset{\mathcal D} {\longrightarrow} N(0,Y),
\end{equation}
where $Y \ge 0$ is a random variable and $N(0,Y)$ is a ``variance mixture'' normal distribution with characteristic function ${\mathbb E} \exp(-Y t^2/2)$, see~\cite{bob_g}. We also note that there is generally no ``loss of mass''  phenomenon for the LIL for trigonometric series with bounded gaps, see~\cite{afu1, afu2}. For further results for trigonometric series with bounded/random  gaps, see~\cite{babera, baza, berbor}.\\

\section{The subsequence principle} \label{sec:subsequ}

The purpose of the previous sections was to illustrate the principle that thin subsequences of the trigonometric system, or thin subsequences of a more general system of dilated functions, exhibit properties which are typical for sequences of independent random variables. However, an analogous principle holds in a much wider framework: it is known that, under suitable technical assumptions, sufficiently thin subsequences of general systems of random variables
behave like genuine independent sequences, in the sense that a general sequence of random variables allows to extract a subsequence showing independent behavior. For example, Gaposhkin
\cite{gap1966,gap1972} and Chatterji~\cite{chaclt,chalilB} proved that if $(X_n)_{n \geq 1}$ is any sequence of random variables
satisfying $\sup_n {\mathbb E} X_n^2<\infty$, then there exist a subsequence
$(X_{n_k})_{k \geq 1}$ and random variables  $X\in L_2$, $Y\in L_1$, $Y \ge 0,$ such that
\begin{equation}\label{CLTmB}
\frac{1}{\sqrt N} \sum_{k\le N} (X_{n_k} - X) \overset{\mathcal D} {\longrightarrow}
N(0,Y)
\end{equation}
and
\begin{equation}\label{LILmB}
\limsup_{N\to\infty} \, \frac{1}{\sqrt{2N \log\log N}} \sum_{k\le N}
(X_{n_k} - X) = Y^{1/2} \qquad \textup{a.s.},
\end{equation}
where as at the end of the previous section $N(0,Y)$ denotes the ``variance mixture''normal distribution with characteristic function
 ${\mathbb E} \exp(-Y t^2/2)$, and where again $\overset{\mathcal D} {\longrightarrow}$ denotes convergence in distribution. A functional (Strassen type) version of (\ref{LILmB})
was proved by Berkes~\cite{BE74B}. By a result of Koml\'os
\cite{koB}, from any sequence $(X_n)_{n \geq 1}$ of random variables satisfying
$\sup_n {\mathbb E} |X_n| < \infty$ one can select a subsequence $(X_{n_k})_{k \geq 1}$
such that
\begin{equation}\label{LLNmB}
\lim_{N\to\infty} \frac{1}{N} \sum_{k\le N} X_{n_k} = X\qquad
\textup{a.s.}
\end{equation}
for some $X\in L_1$.  Chatterji~\cite{chamzB} proved that if $(X_n)_{n \geq 1}$ is a sequence of random variables satisfying $\sup_n {\mathbb E} |X_n|^p < \infty$ for some $0<p<2$, then there exist a subsequence $(X_{n_k})_{k \geq 1}$ and a random variable $X$ with ${\mathbb E} |X|^p<\infty$ such that
\begin{equation}\label{LLNpB}
\lim_{N\to\infty} \frac{1}{N^{1/p}} \sum_{k\le N} (X_{n_k}-X) = 0 \qquad \textup{a.s.}
\end{equation}
These results establish the analogues of the
central limit theorem (CLT), the law of the iterated logarithm (LIL), the
strong law of large numbers (SLLN) and Marczinkiewicz' strong law for subsequences $(X_{n_k})_{k \geq 1}$. Note the
mixed (or randomized) character of  (\ref{CLTmB})--(\ref{LLNpB}):
the limit $X$ in the strong law of large numbers,  the
centering factor $X$ in Marczinkiewicz' strong law,  and the limiting variance $Y$ in the CLT (which
also determines the limsup in the LIL) all become random. For further
limit theorems for subsequences of arbitrary random variable sequences, see
Gaposhkin~\cite{gap1966}. On the basis of these and several other
examples, Chatterji~\cite{chasubB} formulated the following heuristic
principle:

\bigskip\noindent
{\bf Subsequence Principle}. {\sl Let $T$ be a probability limit
theorem valid for all sequences of i.i.d.~random variables belonging
to an integrability class $L$ defined by the finiteness of a norm
$\|\ \cdot \|_L$. Then if $(X_n)_{n \geq 1}$ is an arbitrary (dependent)
sequence of random variables satisfying $\sup_n\|X_n\|_L < + \infty$
then there exists a subsequence $(X_{n_k})_{k \geq 1}$ satisfying $T$ in a
mixed form.}
\bigskip

In a profound study, Aldous~\cite{ald} proved the validity of the subsequence
principle for all distributional and almost sure limit theorems subject to minor technical conditions.
To formulate his results, let $\mathcal M$ denote the
class of probability measures on the Borel sets of $\mathbb R$, equipped with the L\'evy metric. By~\cite{ald}, a subset
$A\subset \mathcal M \times {\mathbb R}^\infty$
is called a {\it limit statute} if:

\bigskip
\begin{enumerate}[(a)]

\item  $P( (\lambda, X_1(\omega), X_2(\omega), \ldots) \in A)=1$ provided $X_1, X_2, \ldots$ are i.i.d.\ random variables with distribution $\lambda$. \\

\item  $(\lambda, x_1, x_2, \ldots )\in A$ and $\sum |x_i-x_i'|<\infty$ implies that $(\lambda, x_1', x_2', \ldots) \in A.$

\end{enumerate}\bigskip

\bigskip

An a.s.\  limit theorem can thus be identified with a limit statute, where the analytic statement of the theorem is expressed by (a), while relation (b) means that a small perturbation of the sequence $X_1, X_2, \ldots$ does not change the validity of the limit theorem. Let us give two examples of limit statutes representing the strong law of large numbers and the law of the iterated logarithm:

\bigskip\noindent
$A_1= \left\{(\lambda, {\bf x})\in A: \lim_{N\to\infty} N^{-1} \sum_{k=1}^N x_k= |\lambda|_1 \right\} \cup \{ (\lambda, {\bf x}): |\lambda|_1=\infty\},$

\bigskip\noindent
$A_2= \left\{(\lambda, {\bf x})\in A:    \limsup_{N\to\infty} (2N\log\log N)^{-1/2} \left(\sum_{k=1}^N x_k-N |\lambda|_1\right) = |\lambda|_2 \right\} \\ \phantom{1111}\cup \{(\lambda, {\bf x}): |\lambda|_2=\infty\}.$

\bigskip\noindent
Here $|\lambda|_1$ and $|\lambda|_2$ denote the mean and variance of $\lambda$ provided they are finite, and we write $|\lambda|_1=\infty$, resp. $|\lambda|_2=\infty$ if $\int_{\mathbb R} |x|d\lambda(x)=\infty$, resp. $\int_{\mathbb R} |x|^2 d\lambda(x)=\infty$.\\

On the other hand, by the definitions in~\cite{ald}, a weak limit theorem for i.i.d.~random variables is a system
$$T =\left(f_1, f_2, \ldots,  \{ G_\lambda, \lambda\in {\mathcal M}_0 \}\right)$$
where

\bigskip
\begin{enumerate}[(a)]
\item ${\mathcal M}_0$ is a measurable subset of $\mathcal M$.\medskip

\item For each $k \ge 1$, $f_k = f_k(\lambda, x_1, x_2, \ldots)$ is a
real function on $\mathcal M \times {\mathbb R}^\infty$, measurable in
the product topology, satisfying the smoothness condition
$$|f_k(\lambda, {\bf x})-f_k(\lambda, {\bf x}')|\le \sum_{k=1}^\infty c_{k, i} |x_i-x_i'|$$
where $0\le c_{k, i}\le 1$ and $\lim_{k\to \infty} c_{k,i}=0$ for each $i$. \medskip

\item For each $\lambda \in {\mathcal M}_0$, $G_\lambda$ is a probability
distribution on the real line such that the map $\lambda \to G_\lambda $
is measurable (with respect to the Borel $\sigma$-field in ${\mathcal M}_0$). \medskip

\item If $\lambda \in {\mathcal M}_0$ and $X_1, X_2,\ldots$ are
independent random variables with common distribution $\lambda $ then
$$ f_k(\lambda,  X_1, X_2, \ldots, ) \overset{\mathcal D} {\longrightarrow} G_\lambda \quad \text{as } k\to\infty.$$
\end{enumerate}
\bigskip

For example, the central limit theorem corresponds to the case when ${\mathcal M}_0$  is the  class of distributions with mean $0$  and finite variance,
\begin{equation} \label{tkcltB}
f_k (\lambda, x_1, x_2, \ldots)= \frac{x_1 + \ldots + x_k - k {\mathbb E}(\lambda)}{
\sqrt k}
\end{equation}
and $G_\lambda= N(0, \text{Var}  (\lambda))$.\\

Let now $(X_n)_{n \geq 1}$ be a sequence of random variables with $\sup_n \|X_n \|_L<\infty$ with any norm $\| \cdot \|_L$  on $\mathbb R$. Then $(X_n)_{n \geq 1}$ is
bounded in probability, i.e.\
$$ \lim_{K\to\infty} P(|X_n|>K)=0 \quad \text{uniformly in } n. $$
By an extension of the Helly--Bray theorem (see e.g.\ \cite{bero}), $(X_n)_{n \geq 1}$ has a subsequence $(X_{n_k})_{k \geq 1}$ having a limit distribution conditionally on any event in the probability space with positive probability, i.e.~for any $A\subset\Omega$
with $P(A) > 0$ there exists a distribution function $F_A$ such that
$$\lim\limits_{k \to\infty} P(X_{n_k} \le  t\mid A) = F_A(t)$$
for all continuity points $t$ of $F_A$. According to the terminology of~\cite{bero}, such a subsequence is called {\it determining}.
Thus when investigating asymptotic properties of sufficiently thin subsequences of sequences $(X_n)_{n \geq 1}$ with bounded norms, we can assume, without loss of generality, that $(X_n)_{n \geq 1}$ itself is determining.
As is shown in~\cite{ald,bero}, for any determining sequence $(X_n)_{n \geq 1}$
there exists a random measure $\mu$ (i.e.\ a measurable
map from the underlying probability space $(\Omega,\mathcal F, \mathcal
P)$ to $\mathcal M$)
such that for any $A$ with $P(A) > 0$ and any
continuity point $t$ of $F_A$ we have
\begin{equation}\label{(4B)} F_A(t) =
{\mathbb E}_A(\mu(-\infty, t]) \end{equation} where ${\mathbb E}_A$ denotes
conditional expectation given $A$. This measure $\mu$ is called the {\it limit
random measure\/} of $(X_n)_{n \geq 1}$; see Section~\ref{sec_new} below for more details.\\

With these preparations, we are now in a position to formulate the subsequence theorems of Aldous.

\begin{thma}[Aldous~\cite{ald}] \label{th8a}
Let $(X_n)_{n \geq 1}$ be a determining sequence with limit random measure $\mu$
and let $A$ be a limit statute. Then there exists a subsequence $(X_{n_k})_{k \geq 1}$
such that for any further subsequence $(X_{m_k})_{k \geq 1} \subset (X_{n_k})_{k \geq 1}$ we have
$$ P( (\lambda(\omega), X_{m_1}(\omega), X_{m_2} (\omega), \ldots) \in A)=1.$$
\end{thma}

\begin{thma}[Aldous~\cite{ald}] \label{th8b}
Let $(X_n)_{n \geq 1}$ be a determining sequence with limit random measure $\mu$
and let
$$T =\left(f_1, f_2, \ldots,  \{ G_\lambda, \lambda\in {\mathcal M}_0 \}\right)$$
be a weak limit theorem. Assume that $P(\mu\in {\mathcal M}_0)=1$. Then there exists a subsequence $(X_{n_k})_{k \geq 1}$
such that for any further subsequence $(X_{m_k})_{k \geq 1} \subset (X_{n_k})_{k \geq 1}$ we have
$$ \lim_{k\to\infty} P( f_k (X_{m_1}(\omega), X_{m_2} (\omega), \ldots \mu(\omega)) \le t)= {\mathbb E} G_{\mu(\omega)} (t)$$
at all continuity points $t$ of the distribution function on the right hand side.
\end{thma}

Writing out Theorem~\ref{th8a} and~\ref{th8b} in the case of the limit statutes $A_1, A_2$ above and the weak limit theorem defined by (\ref{tkcltB}), we get the CLT, LIL and SLLN for thin subsequences of determining sequences, as stated in (\ref{CLTmB}), (\ref{LILmB}), (\ref{LLNmB}) above.\\

The proof of Koml\'os' result (\ref{LLNmB}) exemplifies the technique used in the field of subsequence behavior before Aldous' paper~\cite{ald},  and in particular in proving the results (\ref{CLTmB})--(\ref{LLNpB}) mentioned above.
As Koml\'os showed, if $(X_n)_{n \geq 1}$ is a sequence of random variables with bounded $L_1$ norms, then its sufficiently thin subsequences $(X_{n_k})_{k \geq 1}$ are, after a random centering and small perturbation, an identically distributed martingale difference sequence with finite means and thus, by classical martingale theory, they satisfy the SLLN. Martingale versions of the CLT and LIL yield also relations (\ref{CLTmB}), (\ref{LILmB}) and their functional versions. While this method yields several further limit theorems for lacunary sequences, martingale difference sequences certainly do not satisfy all i.i.d.\ limit theorems in a randomized form and thus the general subsequence principle cannot be proved in such a way. The proof of Theorems~\ref{th8a} and~\ref{th8b} in~\cite{ald}  uses a different way and utilizes near exchangeability properties of subsequences of general sequences of random variables. Let $(X_n)_{n \geq 1}$ be a determining sequence with limit random measure $\mu$ and let $(Y_n)_{n \geq 1}$ be a sequence of random variables, defined on the same probability space as the $X_n$'s, conditionally i.i.d.\ with respect to
$\mu$, with conditional distribution $\mu$. (For the construction of such an $(Y_n)_{n \geq 1}$ one may need to enlarge the probability space.)
Clearly, $(Y_n)_{n \geq 1}$ is exchangeable, i.e.\ for any permutation $\sigma: {\mathbb N} \to {\mathbb N}$ of the positive integers, the sequence $(Y_{\sigma (n)})_{n \geq 1}$ has the same distribution as $(Y_n)_{n \geq 1}$, and it satisfies limit theorems of i.i.d.\ random variables in a mixed form. For example, if ${\mathbb E} Y_1^2<\infty$ and $Y={\mathbb E} (Y_1 \mid \mu)$, $Z=\text{Var}\, (Y_1 \mid \mu)$, then
$$ N^{-1/2} \sum_{k=1}^N (Y_k-Y) \overset{\mathcal D} {\longrightarrow} N(0, Z)$$
and
$$ \limsup_{N\to\infty} \, (2N \log\log N)^{-1/2 } \sum_{k=1}^N (Y_k-Y)=Z^{1/2} \quad \text{a.s.} $$
This principle holds in full generality, i.e.\ for all a.s.\ and distributional
limit theorems in the above formalization. Indeed, if the $Y_n$ are conditionally i.i.d.\
with  respect to $\mu$ and with conditional distribution $\mu$ (a random probability measure on ${\mathbb R}$) and
if $A$ is a limit statute, then
\begin{equation} \label{uncondB}
P( (\mu, Y_1, Y_2, \ldots) \in A  |\mu )(\omega) =P(\mu(\omega), Y_1^*, Y_2^*, \ldots) \in A) \quad \text{a.e.} \emph{}
\end{equation}
where $(Y_n^*)_{n \geq 1}$ is an i.i.d.\ sequence with marginal distribution $\mu (\omega)$. By the definition of limit statute, the last probability in
(\ref{uncondB}) equals 1 and taking expectations we get
$$ P( (\mu, Y_1, Y_2, \ldots) \in A)=1, $$
which is exactly our claim.
Specializing to the case of the limit statutes $A_1$ and $A_2$ above, we get relations (\ref{CLTmB}) and (\ref{LILmB}). A
similar argument works for distributional limit theorems. Now, as is shown in~\cite{ald}, for every $k\ge 1$ we have
\begin{equation}\label {dccB}
(X_{n_1}, X_{n_2}, \ldots X_{n_k}) \overset{\mathcal D} {\longrightarrow} (Y_1, Y_2, \ldots, Y_k) \quad \text{as} \quad  n_1< n_2<\ldots <n_k, \ n_1\to\infty.
\end{equation}
In other words, for large indices the finite dimensional distributions of the sequence $(X_{n_k})_{k \geq 1}$ are close to those of the limiting exchangeable sequence $(Y_k)_{k \geq 1}$ and thus one may expect that limit theorems of $(Y_k)_{k \geq 1}$ (which, as we have just seen, are mixed versions of i.i.d.\ limit theorems) continue to hold for sufficiently thin subsequences $(X_{n_k})_{k \geq 1}$ as well. Of course, a limit theorem for $(X_{n_k})_{k \geq 1}$ can describe a complicated analytic property of the infinite vector $(X_{n_1}, X_{n_2}, \ldots, X_{n_k}, \ldots)$ which does not follow from the weak convergence of the finite dimensional distributions of the sequence, but with a suitable thinning procedure and delicate analytic arguments, Aldous showed an infinite dimensional extension of (\ref{dccB}), leading to the validity of Theorems~\ref{th8a} and~\ref{th8b}.\\

Although the theorems of Aldous are of exceptional generality, there are important results for lacunary sequences which are not covered by them.  As was shown by Gaposhkin~\cite{gap1966}, for every uniformly bounded sequence $(X_n)_{n \geq 1}$ of  random variables there exists a subsequence $(X_{n_k})_{k \geq 1}$ and bounded random variables $X$ and $Y\ge 0$ such that for any numerical sequence $(a_n)_{n \geq 1}$  satisfying
\begin{equation}\label{coeffB}
A_N:=\sum_{k=1}^N a_k^2\to\infty, \quad a_N=o(A_N^{1/2})
\end{equation}
we have
\begin{equation}\label{cltakB}
\frac{1}{A_N} \sum_{k=1}^N a_k (X_{n_k}-X) \overset{\mathcal D} {\longrightarrow} N(0, Y),
\end{equation}
and if the second relation of (\ref{coeffB}) is replaced  by
\begin{equation}\label{lilakB}
a_N=o(A_N/(\log\log A_N)^{1/2})
\end{equation}
then we have
\begin{equation}\label{LILkolmB}
\limsup_{N\to\infty} \frac{1}{\sqrt{2A_N^2 \log\log A_N}} \sum_{k=1}^N a_k (X_{n_k}-X) =Y^{1/2} \quad \text{a.s.}
\end{equation}
The difference of these results from (\ref{CLTmB}) and (\ref{LILmB}) is that
in the CLT and LIL we have weighted sums $\sum_{k=1}^N a_k (X_{n_k}-X)$ instead of ordinary sums $\sum_{k=1}^N (X_{n_k}-X)$. For every fixed coefficient sequence $(a_n)_{n \geq 1}$ the CLT and LIL in (\ref{cltakB}) and (\ref{lilakB}) follow from Theorems~\ref{th8a} and~\ref{th8b}, but the subsequence $(X_{n_k})_{k \geq 1}$ provided by the proofs  depends on $(a_k)_{k \geq 1}$ and it is not clear that we can select a subsequence $(X_{n_k})_{k \geq 1}$ satisfying (\ref{cltakB}) and (\ref{lilakB}) simultaneously for all considered coefficient sequences $(a_k)_{k \geq 1}$.\\

Another important situation not covered by Aldous' general theorems is when we investigate permutation-invariance of limit theorems for subsequences. Since the asymptotic properties of an exchangeable sequence $(Y_n)_{n \geq 1}$ do not change after any permutation of its terms, it is natural to expect that the conclusions in Theorem~\ref{th8a} and~\ref{th8b} remain valid after an arbitrary permutation of the subsequence $(X_{n_k})_{k \geq 1}$ in the theorems. However, the proofs of Theorem~\ref{th8a} and~\ref{th8b} are not permutation-invariant and it does not follow that, e.g.,  any sequence $(X_n)_{n \geq 1}$ of random variables with bounded $L_1$ norms contains a subsequence $(X_{n_k})_{k \geq 1}$ satisfying the strong law of large numbers  after any permutation of its terms. Using ad hoc methods, the latter result has been proved by Berkes~\cite{be1985B} and another classical case, namely the unconditional a.e.\ convergence of series $\sum c_k (X_{n_k}-X)$ under $\sum c_k^2<\infty$ for subsequences $(X_{n_k})_{k \geq 1}$ of $L_2$ bounded sequences $(X_n)_{n \geq 1}$, has been settled by Koml\'os~\cite{kopermB} (see~\cite{ald} for another  proof via exchangeability). It clearly would be desirable to provide further general results in this direction.\\

We now formulate some structure theorems for lacunary sequences enabling one to handle problems of the kind discussed above.  Recall that if $(X_n)_{n \geq 1}$ is a determining sequence with limit random measure $\mu$  and $(Y_n)_{n \geq 1}$ is a sequence conditionally i.i.d.\ with respect to the $\sigma$-algebra generated by $\mu$ and with conditional marginal distributions $\mu$, then there exists a subsequence $(X_{n_k})_{k \geq 1}$ such that (\ref{dccB}) holds. This shows that, in some sense, for large indices the sequence $(X_{n_k})_{k \geq 1}$ resembles the sequence $(Y_k)_{k \geq 1}$, but this property is far too weak to deduce limit theorems for $(X_{n_k})_{k \geq 1}$ from those valid for the exchangeable sequence $(Y_k)_{k \geq 1}$.  The following theorem, proved by Berkes and P\'eter~\cite{bp},  shows that with a suitable choice of the subsequence $(n_k)_{k \geq 1}$, the variables $(X_{n_k})_{k \geq 1}$ can be chosen to be close to the $Y_k$ in a pointwise sense. We call a sequence $(X_n)_{n \geq 1}$ of random variables $\e$-{\it exchangeable} if on the same probability space  there exists an exchangeable sequence $(Y_n)_{n \geq 1}$ such that $P(|X_n-Y_n|\ge \e)\le \e$ for all $n$. Then we have

\begin{thma}[Berkes and P\'eter~\cite{bp}] \label{th8c}
Let $(X_n)_{n \geq 1}$ be a sequence of random variables bounded in probability, and let $(\ve_n)_{n \geq 1}$ be a sequence of positive reals tending to zero. Then, if the underlying probability space is large enough, thee exists a subsequence $(X_{n_k})_{k \geq 1}$ such that, for all $l\ge 1$, the sequence $X_{n_l}, X_{n_{l+1}}, \ldots$ is $\ve_l$-exchangeable.
\end{thma}

Note that Theorem~\ref{th8c} provides a different approximating exchangeable sequence $(Y_j^{(l)})_{j \geq 1}$ for each tail sequence $(X_{n_l}, X_{n_{l+1}}, \ldots)$, with termwise approximating error $\ve_l$. The following theorem describes precisely the structure of the  the sequences $(Y_j^{(l)})_{j\geq 1}$.

\begin{thma}[Berkes and P\'eter~\cite{bp}] \label{th8d}
Let $(X_n)_{n \geq 1}$ be a determining sequence of random variables, and let
$(\ve_n)_{n \geq 1}$ be a sequence of positive reals. Then there exists a
subsequence $(X_{m_k})_{k \geq 1}$ and a sequence $(Y_k)_{k \geq 1}$ of discrete random variables
such that
\begin{equation}\label{2.1B}
P\bigl( |X_{m_k} - Y_k| \geq \ve_k \bigr) \leq \ve_k \quad k = 1,2
\dots, \end{equation} and for each $k > 1$ the atoms of the finite
$\sigma$-field $\sigma \{ Y_1, \dots, Y_{k - 1} \}$ can be divided
into two classes $\Gamma_1$ and $\Gamma_2$ such that the following holds. Firstly,
\begin{equation}\label{rmiB}
\sum_{B \in \Gamma_1} P(B) \leq \ve_k.
\end{equation}
Secondly, for any $B \in \Gamma_2$ there exist
$P_B$-independent random variables $\{ Z^{(B)}_j, j = k, k + 1, \dots \}$
defined on $B$ with common distribution function $F_B$ such that
\begin{equation}\label{2.2B}
P_B \bigl( |Y_j - Z^{(B)}_j | \geq \ve_k \bigr) \leq \ve_k, \qquad j
= k, k + 1, \dots
\end{equation}
Here $F_B$ denotes the limit distribution of $(X_n)_{n \geq 1}$ relative to
$B$ and $P_B$ denotes conditional probability given~$B$.
\end{thma}

We now give applications of Theorem~\ref{th8d} to the problems discussed above.
First we note that using Theorem~\ref{th8d} it is a simple exercise to prove, for suitable subsequences of a uniformly bounded sequences
$(X_n)_{n \geq 1}$, the weighted CLT and LIL in (\ref{cltakB}), (\ref{LILkolmB}) simultaneously for all permitted coefficient
sequences $(a_n)_{n \geq 1}$. Next we give a permutation-invariant form of Theorem~\ref{th8b} for distributional limit theorems.

\begin{defn}
We call the weak limit theorem $T = (f_1,
f_2, \ldots, S, \, \{ G_\mu, \mu\in {\mathcal M}_0\})$ regular if there exist
sequences $p_k \le q_k$ of positive integers tending to $+\infty$
and a sequence $\omega_k\to +\infty$ such that

\begin{enumerate}[(i)]
\item $f_k(\lambda,  x_1, x_2, \ldots)$ depends only on
$\lambda, x_{p_k}, \ldots, x_{q_k}$. \medskip

\item $f_k$ satisfies the Lipschitz condition

$$|f_k(\lambda, x_{p_k}, \ldots, x_{q_k}) - f_k (\lambda', x'_{p_k},\ldots,
x'_{q_k})| \le$$
$$\le \frac{1}{\omega_k} \sum^{q_k}_{i = p_k} |x_i - x'_i|^\alpha
+ \varrho^* (\lambda,\lambda')$$ for some $0 < \alpha \le 1$, where
$\varrho^*$ is a metric on ${\mathcal M}_0$ generating the same topology as the
Prohorov metric $\varrho$.
\end{enumerate}
\end{defn}

For example, the central limit theorem can be formalized by the functions
$$
f_k(\lambda, x_{[k^{1/4}]}, \ldots, x_k) = \frac{x_{[k^{1/4} ]} +\ldots +
x_k - k {\mathbb E} (\lambda)}{\sqrt k},
$$
leading to a regular limit theorem. Note that originally we formalized the CLT with the functions $f_k$ in
(\ref{tkcltB}) containing all variables $x_1, x_2, \ldots$, but
under bounded second moments the first $k^{1/4}$ terms here are irrelevant and hence we can always switch
to the present version. The same procedure applies in the general case.

\begin{thma}[Aistleitner, Berkes and Tichy~\cite{abt_perm2}] \label{th8e}
Let $(X_n)_{n \geq 1}$ be a determining sequence
with limit random measure $\tilde{\mu}$. Let $T = (f_1, f_2, \ldots,
S, \{ G_\mu, \mu\in {\mathcal M}_0\} )$ be a regular weak limit theorem and
assume that $P(\tilde{\mu} \in {\mathcal M}_0) = 1$. Then there exists a
subsequence $(X_{n_k})_{k \geq 1}$ such that for any permutation $(X^*_k)_{k \geq 1}$
of $(X_{n_k})_{k \geq 1}$ we have
\begin{equation} \label{pB}
f_k(X_1^*,  X_2^*, \ldots, \tilde{\mu} ) \to_d \int G_{\tilde{\mu}} dP
.
\end{equation}
\end{thma}

In case of the CLT formalized above, assuming $\sup_n {\mathbb E}
X_n^2<+\infty$ implies easily that $\tilde{\mu}$ has finite variance
almost surely, and thus denoting its mean and variance by $X$ and
$Y$, respectively, we see that the integral in (\ref{pB}) is the
distribution $N(0, Y)$. Hence (\ref{pB}) states in the present case
that
$$\frac{1}{\sqrt{N}} \sum_{k=1}^N (X_k^*-X) \overset{\mathcal D} {\longrightarrow}  N(0, Y),$$
which is the permutation-invariant form of the CLT.  \\

Concerning a.s.\ limit theorems, a permutation-invariant form of the strong law of large numbers
for subsequences of an $L_1$-bounded sequence was proved, as already mentioned, in Berkes~\cite{be1985B},
and a similar argument yields the corresponding result for the LIL.
No permutation-invariant version of the general result in Theorem~\ref{th8a} has been proved in the literature, but there is
no need for that, since a.s.\ limit theorems can be reformulated in a distributional
form and thus the proof of Theorem~\ref{th8b} applies with obvious changes. For illustration, we give here the
reformulation of the LIL:

\begin{thma}
Let $(X_n)_{n \geq 1}$ be a sequence of random variables with ${\mathbb E} |X_n| \le 1$,
$n=1,2,\ldots$ Put $S_n = \sum_{i=1}^n X_i$, $S_{k,
l}=\sum_{i=k+1}^l X_i$, and $L(N)=(2N\log\log N)^{1/2}.$  Then $\limsup_{N\to\infty} S_N/L(N)=1$
a.s.\ iff for any $\ve>0$ there exists a sequence $m_1<m_2<\cdots\ $
of positive integers such that $m_k\ge 5^k$ and
$$
P \left( \max_{m_k\le j\le m_{k+1}} \frac{S_{k, j}}{L(j)} >1+\ve \right) \le 2^{-k}, \qquad k\ge k_0,
$$
and
$$
P \left( \max_{m_k\le j\le m_{k+1}} \frac{S_{k, j}}{L(j)} <1-\ve \right) \le 2^{-k}, \qquad
k\ge k_0.
$$ 
\end{thma}

It is worth pointing out that given a sequence $(X_n)_{n \geq 1}$ of random variables, finding a subsequence $(X_{n_k})_{k \geq 1}$ satisfying the permutation-invariant form of some limit theorem generally requires a much faster growing sequence $(n_k)_{k \geq 1}$ than to find a subsequence to satisfy the original limit theorem. This is a phenomenon which also occurs for lacunary trigonometric sums or lacunary sums of dilated functions; compare the last paragraph of Section~\ref{sec:CLT} above.\\

In conclusion we note that if $(X_n)_{n \geq 1}$  is a sequence of random variables with finite means over the probability space $(0, 1)$ equipped with the Borel $\sigma$-algebra and Lebesgue measure such that for all $n\ge 1$ and $(a_1, \ldots, a_n)\in {\mathbb R}^n$ we have
\begin{equation} \label{lp}
C_1 \left(\sum_{k=1}^n |a_k|^p\right)^{1/p} \le {\mathbb E} \left| \sum_{k=1}^n a_k X_k \right| \le C_1 \left(\sum_{k=1}^n |a_k|^p\right)^{1/p}
\end{equation}
for some $p\ge 1$ and positive constants $C_1, C_2$, then the closed subspace of $L_1(0, 1)$ spanned by the $X_n$ is isomorphic with the $\ell_p$ space (Hilbert space if $p=2$). Relation (\ref{lp}) holds, in particular, if the $X_n$ are i.i.d.\ symmetric $p$-stable random variables with $p>1$, i.e.\ their characteristic function (Fourier transform) is given by $\exp(-c|t|^p)$ with some $c>0$.
Thus applying the subsequence principle to the ``limit theorem'' (\ref{lp}) provides important information on the subspace structure of $L_1(0, 1)$. Using this method, Aldous~\cite{ald81B} proved the famous conjecture that every infinite dimensional closed subspace of $L_1(0, 1)$ contains an isomorphic copy of $\ell_p$  for some $1\le p\le 2$. For a further application of this method, see an improvement of the classical theorem of Kadec and Pe{\l}czy\'nski~\cite{kape} on the subspace structure on $L_p$, $p>2$, in Berkes and Tichy~\cite{bt}.

\section{New results: Exact criteria for the central limit theorem for subsequences} \label{sec_new}

By the classical resonance theorem of Landau~\cite{lan}, for a real sequence $(x_n)_{n \geq 1}$ the series $\sum_{n=1}^\infty a_n x_n$
converges for all $(a_n)_{n \geq 1} \in \ell_p$ $(1\le p\le \infty)$ if and only if $(x_n)_{n \geq 1} \in \ell_q$, where  $1/p+1/q=1$. A deep extension of this result
to the case of function series was given by Nikishin~\cite{nik}.
We call a sequence $(f_n)_{n \geq 1}$ of measurable functions on $(0,1)$ a {\it convergence system in measure for $\ell_p$} if for any real sequence
$(a_n)_{n \geq 1} \in \ell_p$ the series  $\sum_{n=1}^\infty a_n f_n$ converges in measure.
In the case $p=2$ Nikishin proved the following result.

\begin{thma}[Nikishin~\cite{nik,nik2}] \label{tha}
A function system $(f_n)_{n \geq 1}$ over $(0,1)$ is a convergence system in measure for $\ell_2$ if and only if for any
$\ve>0$ there exists a measurable set $A_\ve \subset (0, 1)$ with measure  exceeding $1-\ve$ and a constant $K_\ve>0$ such that for all $N\ge 1$ and all $(a_1, \ldots, a_N) \in {\mathbb R}^N$ we have
\begin{equation} \label{ncond}
\int_{A_\ve} \left(\sum_{n=1}^N a_n f_n\right)^2 dx \le K_\ve \sum_{n=1}^N a_n^2.
\end{equation}
\end{thma}

The sufficiency of \eqref{ncond} is obvious, so the essential (and highly remarkable) statement is the converse: if a sequence $(f_n)_{n \geq 1}$ is a convergence system in measure for $\ell_2$, then, except for a subset of $(0,1)$ with arbitrary small measure, $(f_n)_{n \geq 1}$ behaves like an orthonormal sequence.\\

In the previous section we discussed the subsequence principle stating that sufficiently thin subsequences of arbitrary sequences of random variables, subject to mild boundedness conditions, satisfy ``all'' limit theorems for i.i.d.\ random variables in a mixed (randomized) form. A typical special case of this principle is the following result:

\begin{thma}[Gaposhkin~\cite{gap1972}] \label{thb}
Let $(X_n)_{n \geq 1}$ be a sequence of random variables satisfying
\begin{equation}\label{1.1}
\sup_n {\mathbb E} X^2_n < + \infty .
\end{equation}
Then there exists a subsequence $(X_{n_k})_{k \geq 1}$ together with random variables $X$ and $Y \geq 0$ such that for any further subsequence $(X_{m_k})_{k \geq 1}$ of
$(X_{n_k})_{k \geq 1}$ we have
\begin{equation}\label{gapclt}
\frac{1}{\sqrt N} \sum_{k=1}^N (X_{m_k} - X) \overset{\mathcal D} {\longrightarrow} N(0, Y),
\end{equation}
where $N(0,Y)$ denotes the ``variance mixture'' normal distribution with characteristic
function ${\mathbb E} \exp(-Y t^2/2)$.
\end{thma}

If $(X_n^2)_{n \geq 1}$ is uniformly integrable then by well-known compactness results (see e.g.\ \cite{ds}) there
exist a subsequence $(X_{m_k})_{k \geq 1}$ and random variables $X\in L_2$
and $Y\in L_{1/2}$, $Y\ge 0$, such that
\begin{equation}\label{comp}
X_{m_k}\to X \ \text{weakly in}  \ L_2, \qquad (X_{m_k}-X)^2 \to Y^2 \ \text{weakly in} \ L_1.
\footnote{A sequence $(\xi_n)_{n \geq 1}$ of random variables in $L_p$, $p\ge 1$, is said to converge weakly to $\xi\in L_p$
if ${\mathbb E} (\xi_n\eta)\to {\mathbb E} (\xi\eta)$ for any $\eta\in L_q$, where $1/p+1/q=1$. This type of convergence should not be confused
with weak convergence of probability measures and distributions, called generally convergence in distribution,
and denoted by $\overset{\mathcal D} {\longrightarrow}$.}
\end{equation}
As Gaposhkin~\cite{gap1966} showed, in this case the random variables $X, Y$ in \eqref{gapclt} can be chosen as in \eqref{comp}.\\

In Theorem~\ref{thb}, condition (\ref{1.1}) is not necessary: simple examples show (see below) that there exist sequences $(X_n)_{n\ge 1}$ of random variables without any finite moments, but having subsequences satisfying (\ref{gapclt}). The purpose of this section is to give necessary and sufficient conditions for the existence of subsequences $(X_{n_k})_{k \ge 1}$ satisfying the randomized CLT (\ref{gapclt}), and it will turn out that our conditions have the same character as Nikishin's conditions for the existence of a subsequence being a convergence system, i.e.\ ``nice'' behavior of the sequence on subsets of the probability space with measure as close to 1  as we wish.\\

To formulate our results, call a sequence $(X_n)_{n \geq 1}$ of random variables {\it nontrivial} if it has no subsequence converging with positive probability. It is easily seen that for non-degenerate sequences the random variable $Y$ in Theorem~\ref{thb} is almost surely positive and Gaposhkin's theorem can be rewritten in a form involving a pure (i.e.\ not mixed) Gaussian limit distribution.

\addtocounter{thma}{1}

\begin{thma} \label{thc}
Let $(X_n)_{n \geq 1}$ be a nontrivial sequence of random variables satisfying \eqref{1.1}.
Then there exists a subsequence $(X_{n_k})_{k \geq 1}$ and random variables $X, Y$
with $Y > 0$ such that for all subsequences $(X_{m_k})_{k \geq 1}$ of
$(X_{n_k})_{k \geq 1}$ and for any set $A$ in the probability space with
$P(A) > 0$ we have
\begin{equation}\label{1.2}
P_A \left( \frac{\sum_{k=1}^N (X_{m_k} - X)}{Y \sqrt
N} < t \right) \to \Phi (t) \quad \text{for all } t.
\end{equation}
Here $P_A$ denotes the conditional probability with respect to~$A$, and $\Phi$ is the cumulative distribution function of the standard normal distribution.
\end{thma}

The nontriviality of $(X_n)_{n \geq 1}$ is assumed here to avoid degenerate
cases. If $X_{n_k} \to X$ on some set $A$ with positive
probability then for any sufficiently thin subsequence
$(X_{m_k})_{k \geq 1}$ of $(X_{n_k})_{k \geq 1}$ we have $\sum |X_{m_k} - X| < +\infty$
a.s.\ on $A$, and consequently
$$
a^{-1}_N \sum_{k=1}^N (X_{m_k} - X) \to 0 \quad \text{a.s. on
}\ A
$$
for any norming sequence $a_N \to \infty$ (random or not). Since
for any sequence $(X_n)_{n \geq 1}$ satisfying \eqref{1.1} (and in fact any tight
sequence $(X_n)_{n \geq 1}$) there is a subsequence $(X_{n_k})_{k \geq 1}$ and a
measurable partition $A \cup B$ of the probability space such
that $X_{n_k}$ converges on $A$ and is nontrivial on $B$, there
is no loss of generality in assuming that $(X_n)_{n \geq 1}$ is nontrivial.\\

Clearly, if $(X_n)_{n \geq 1}$ satisfies the conclusion of Theorem~\ref{thc}, then so does
the sequence $(X_n+2^{-n}Z)_{n\ge 1}$ for any a.s.\ finite random variable Z, and thus the assumption
\eqref{1.1} is, as stated above, not necessary in Theorem~\ref{thc}. Below we will give necessary and sufficient condition for the CLT for lacunary subsequences
of a given sequence $(X_n)_{n\ge 1}$ of random variables without any moment assumption on $(X_n)_{n \geq 1}$.
To formulate our results, let us note that if all subsequences
$(X_{m_k})_{k \geq 1}$ of a sequence $(X_n)_{n \geq 1}$ satisfy \eqref{1.2} for some random variables
$X,Y$, then $(X_n)_{n \geq 1}$ is bounded in probability (see Lemma 2 below).
As mentioned in the previous section,
every sequence $(X_n)_{n \geq 1}$ of random variables bounded in probability has a
subsequence $(X_{n_k})_{k \geq 1}$ which has a limit distribution relative to
every set $A$ of the probability space with $P(A) > 0$. Such a sequence was called {\it determining}.
This concept is the same as that of stable convergence, introduced by
R\'enyi~\cite{renyi}; our terminology follows that of functional analysis.
Hence in our investigations we can assume without loss of generality
that the original sequence $(X_n)_{n \geq 1}$ is determining. Now if $(X_n)_{n \geq 1}$ is
determining and $F_A$ denotes its limit distribution relative to the set $A$,
then, as we noted in the previous section, there exists a random measure
$\mu$ (called the limit random measure of $(X_n)$) such that
\begin{equation}\label{1.4}
F_A (t) = {\mathbb E}_A (\mu (-\infty, t])
\end{equation}
for any continuity point $t$ of $F_\Omega$, where ${\mathbb E}_A$ denotes
conditional expectation relative to~$A$. Let $F_\bullet$ denote the distribution function of
$\mu$; we shall call it the {\it limit random distribution\/} of~$(X_n)_{n \geq 1}$. We can state now our first new theorem.

\begin{thm} \label{th1}
Let $(X_n)_{n \geq 1}$ be a nontrivial sequence of random variables.
Then the following statements are equivalent:\\

\begin{enumerate}[A)]
\item There exist a subsequence $(X_{n_k})_{k \geq 1}$ and random variables $X, Y$
with $Y > 0$ such that \eqref{1.2} holds for all subsequences
$(X_{m_k})_{k \geq 1}$ of $(X_{n_k})_{k \geq 1}$ and for any set $A \subset \Omega$
with $P(A) > 0$.\\

\item For every $\ve > 0$ there is a subsequence $(X_{n_k})_{k \geq 1}$ and a
set $A_\varepsilon \subset \Omega$ with $P(A_\varepsilon) \geq 1 - \ve$ such that
\begin{equation}\label{1.5}
\sup_k \int\limits_{A_\varepsilon} X^2_{n_k} d P < + \infty .
\end{equation}
\end{enumerate}

\medskip
If $(X_n)_{n \geq 1}$ is determining, then two further equivalent statements are:\\

\begin{enumerate}[A)]    \addtocounter{enumi}{2}

\item We have
\begin{equation} \label{rmC}
\int\limits^{+\infty}_{-\infty} x^2 dF_\bullet(x) < + \infty \quad \text{almost surely.}
\end{equation}

\item For every $\ve > 0$ there exists a set $A_\varepsilon \subset \Omega$
with $P(A_\varepsilon) \geq 1 - \ve$ such that
\begin{equation}\label{1.6}
\int\limits^{+\infty}_{-\infty} x^2 dF_{A_\varepsilon}(x) < + \infty.
\end{equation}

\end{enumerate}
\end{thm}

Our second new theorem characterizes sequences $(X_n)_{n \geq 1}$ for which \eqref{1.2}
holds with $X \in L_2$, $Y \in L_{1/2}$.

\begin{thm} \label{th2}
Let $(X_n)_{n \geq 1}$ be a nontrivial sequence of random variables defined on an
atomless probability space $(\Omega, \mathcal F, P)$. Then the
following statements are equivalent:

\begin{enumerate}[A)]
\item There exists a subsequence $(X_{n_k})_{k \geq 1}$ and random variables $X,Y$ with
$Y > 0$, $X \in L_2$, $Y \in L_{1/2}$ such that \eqref{1.2} holds for
all subsequences $(X_{m_k})_{k \geq 1}$ of $(X_{n_k})_{k \geq 1}$ and all sets $A
\subset \Omega$ with $P(A) > 0$.\\

\item There exists a subsequence $(X_{n_k})_{k \geq 1}$ and sequences
$(Y_k)_{k \geq 1}$, $(\tau_k)_{k \geq 1}$ of random variables satisfying
\begin{equation}\label{1.7}
X_{n_k} = Y_k + \tau_k,
\end{equation}
where
\begin{equation}\label{1.8}
\sup_k {\mathbb E} Y^2_k < + \infty, \quad \sum_k |\tau_k| < + \infty \qquad
\text{a.s.}
\end{equation}
\end{enumerate}

\medskip
If $(X_n)_{n \geq 1}$ has a limit distribution $F$, then a third equivalent
statement is:\\

\begin{enumerate}[A)]    \addtocounter{enumi}{2}

\item We have
\begin{equation}\label{rmCC}
\int\limits^{+\infty}_{-\infty} x^2 dF(x) < + \infty.
\end{equation}
\end{enumerate}
\end{thm}

In other words, for the validity of \eqref{1.2} with $X \in L_2$, $Y
\in L_{1/2}$, assumption \eqref{1.1} is necessary and sufficient after a small
perturbation of $(X_n)_{n \geq 1}$, and for identically distributed $(X_n)_{n \geq 1}$
even this perturbation is not needed. A particularly simple case when $X \in
L_2$, $Y \in L_{1/2}$ is satisfied is when $X,Y$ are nonrandom.\\

A trivial example showing the difference between condition (D)
of Theorem~\ref{th1} and condition (C) of Theorem~\ref{th2} is the following.
Let $\{ H_k, k \geq 1\}$ be a partition of the probability space
with $P(H_k) = 2^{-k}$ for $k = 1,2, \dots$, and let $(X_n)_{n \geq 1}$ be a
sequence of random variables on this space which is conditionally i.i.d.\
given each $H_k$ with mean $0$ and variance $2^k$. Then $(X_n)_{n \geq 1}$
is nontrivial, determining and clearly satisfies condition (D)
of Theorem~\ref{th1}, but since it is identically distributed (in fact
exchangeable) and since ${\mathbb E} X^2_1 = +\infty$, condition (C) of Theorem~2
is not satisfied.

\subsection{Some lemmas}

The key for the proof of our theorems is a general structure
theorem for lacunary sequences which was proved in~\cite{bp}, and which was stated as Theorem~\ref{th8d} in the previous section. Furthermore, we need the following lemmas.

\begin{lemma} \label{lemma2}
Let $(X_n)_{n \geq 1}$ be a sequence of random variables such that for
some random variables $X, Y$ with $Y > 0$ and for all subsequences
$(X_{n_k})_{k \geq 1}$ we have
\begin{equation}\label{2.3}
\frac{\sum\limits_{k=1}^N (X_{n_k} - X)}{Y \sqrt N}
\overset{\mathcal D} \longrightarrow N(0,1).
\end{equation}
Then $(X_n)_{n \geq 1}$ is bounded in probability.
\end{lemma}

\begin{proof}
Clearly \eqref{2.3} implies that the sequence $(X_{n_N} - X) / (Y
\sqrt{N})$ is bounded in probability as $N \to \infty$, and thus $X_{n_N} / \sqrt N$
is bounded in probability for any subsequence $(n_k)_{k \geq 1}$. If
$(X_n)_{n \geq 1}$ were not bounded in probability then one could find a
subsequence $(m_k)_{k \geq 1}$ and a constant $c > 0$ such that
$P(|X_{m_k}| \geq k) \geq c$ for $k = 1,2, \dots$, i.e.\
$X_{m_k} / \sqrt k$ would not be bounded in probability, a
contradiction.
\end{proof}

\begin{lemma} \label{lemma3}
Let $(X_n)_{n \geq 1}$ be a sequence of random variables and assume that for some
random variables $X$ and $Y > 0$ and all sets $A \subset \Omega$ with
$P(A) > 0$ we have
\begin{equation}\label{2.4}
P_A \left( \frac{\sum\limits_{k=1}^N (X_k - X)}{Y \sqrt N} <
t\right) \to \Phi (t) \quad \text{for all } t.
\end{equation}
Assume further that \eqref{2.4} remains valid if we replace $X, Y$ by
some random variables $X^*$ and $Y^* > 0$. Then $X = X^*$ a.s.\ and $Y =
Y^*$ a.s.
\end{lemma}

\begin{proof}
From the assumption it follows that the sequences
$$
N^{-1/2}
\sum\limits_{k=1}^{N} (X_k - X) \qquad \text{and} \qquad N^{-1/2} \sum_{k=1}^{N}
(X_k - X^*)
$$
are bounded in probability, and thus the same holds
for their difference $\sqrt N(X - X^*)$, whence $X = X^*$ a.s. To
prove $Y = Y^*$, fix $c > 1$ and set $A = \{ Y^* \geq cY\}$. If
$P(A) > 0$ then clearly we cannot have both \eqref{2.4} and the
analogous relation with $Y$ replaced by $Y^*$. Thus $P(A) > 0$
for all $c > 1$ whence $Y^* \leq Y$ a.s. The same argument
yields $Y \leq Y^*$ a.s., completing the proof.
\end{proof}

\begin{lemma} \label{lemma4}
Let $X_1, X_2, \dots, X_n$ be i.i.d.\ random variables with
distribution function $F$ and set $S_{n} = X_1 + \cdots +
X_n$. Then for any $t > 0$ we have
\begin{equation}\label{2.5}
P(|S_n| \leq 2t) \leq A \frac{t}{\sqrt n} \left( \int\limits_{|x| \leq t} x^2 dF(x) - 2 \left(
\int\limits_{|x| \leq t} x dF(x)\right)^2 \right)^{-1/2},
\end{equation}
provided the difference on the right-hand side is positive and
$\int\limits_{|x| \leq t} dF(x) \geq 1/2$. Here $A$ is an absolute constant.
\end{lemma}

\begin{proof}
Let $F^*$ denote the distribution function obtained
from $F$ by symmetrization. From a well-known concentration
function inequality of Esseen~\cite[Theorem 2]{es} it follows
that the left-hand side of \eqref{2.5} cannot exceed
$$
A_1 \frac{t}{\sqrt n} \left( \int\limits_{|x| \leq 2t} x^2 dF^* (x)
\right)^{-1/2},
$$
where $A_1$ is an absolute constant. Hence to prove \eqref{2.5} it
suffices to show that $\int\limits_{|x| \leq t} dF(x) \geq 1/2$ implies
\begin{equation}\label{2.6}
\int\limits_{|x| \leq 2t} x^2 dF^*(x) \geq \int\limits_{|x| \leq t} x^2 dF(x)
- 2 \left(\int\limits_{|x| \leq t} x dF(x) \right)^2.
\end{equation}
Let $\xi$ and $\eta$ be independent random variables with distribution
function $F$, and set
$$
C = \{ | \xi - \eta | \leq 2t \}, \qquad D = \{ |\xi | \leq t,
|\eta| \leq t \}.
$$
Then
\begin{align*}
\int\limits_{|x| \leq 2t} x^2 dF^*(x)
&= \int\limits_C (\xi - \eta)^2 dP \\
& \geq \int\limits_D (\xi - \eta)^2 dP \\
&= 2 \int\limits_{|\xi| \leq t} \xi^2 dP \cdot P(|\eta | \leq t) - 2
\left( \int\limits_{|\xi| \leq t} \xi dP \right)^2 \\
&\geq \int\limits_{|\xi|\leq t} \xi^2 dP - 2 \left( \int\limits_{|\xi| \leq
t} \xi dP\right)^2,
\end{align*}
provided $P(|\eta| \leq t) \geq 1/2$. Thus \eqref{2.6} is valid.
\end{proof}

\begin{lemma} \label{lemma5}
Let $(\Omega, \mathcal F, P)$ be an atomless probability space and
$X_1, X_2, \dots$ a sequence of random variables on $(\Omega, \mathcal F, P)$
with limit distribution $F$. Then there exist a subsequence
$(X_{n_k})_{k \geq 1}$ and sequences $(Y_k)_{k \geq 1}$ and $(\tau_k)_{k \geq 1}$ of random variables on
$(\Omega, \mathcal F, P)$ such that $X_{n_k} = Y_k + \tau_k$, $k =
1,2, \dots$, such that the random variables $Y_k$ have distribution function $F$,
and such that $\sum_k |\tau_k| < + \infty$ a.s.
\end{lemma}

\begin{proof}
Let $(\hat X_n)_{n \geq 1}$ be discrete random variables such that $P(|X_n - X'_n| \geq
2^{-n}) \leq 2^{-n}$,  $n = 1,2, \dots$, and denote by $F_n$ the
distribution function of $X_n$. Clearly $F_n \to F$ and thus
$\ve_n : = \varrho (F_n, F) \to 0$, where $\varrho$ denotes the
Prohorov distance. By a theorem of Strassen [8] there exists a
probability measure $\mu_n$ on $\mathbb{R}^2$ with marginals $F_n$ and
$F$ such that
$$
\mu_n((x,y) : |x - y| \geq \ve_n) \leq \ve_n.
$$
Let $c$ be a possible value of $\hat X_n$. Since the probability
space restricted to $A = \{ \hat X_n = c\}$ is atomless, there
exists a random variable $V_n$ on this space such that
$$
P_A(V_n < t) = \frac{\mu_n((x,y) : x = c, y < t)}{\mu_n((x,y) : x = c)}
$$
for all $t$. Carrying out this construction for all possible
values of $c$ in the range of $\hat X_n$, we get a random variable $V_n$
defined on the whole probability space such that the joint
distribution of $\hat X_n$ and $V_n$ is $\mu_n$. Clearly the
distribution of $V_n$ is $F$ and $P(|\hat X_n - V_n| \geq \ve_n)
\leq \ve_n$. Choosing $(n_k)_{k \geq 1}$ so that $\ve_{n_k} \leq 2^{-k}$ we
get
$$
P\bigl( |\hat X_{n_k} - V_{n_k} | \geq 2^{-k} \bigr) \leq
2^{-k},
$$
i.e.
$$
P \bigl( |X_{n_k} - V_{n_k} | \geq 2 \cdot 2^{-k} \bigr) \leq 2
\cdot 2^{-k}
$$
and thus $\sum_k | X_{n_k} - V_{n_k}| < + \infty$ a.s.\ by the
Borel--Cantelli lemma. Thus the decomposition $X_{n_k} = Y_k +
\tau_k$, where $Y_k = V_{n_k}$ and $\tau_k = X_{n_k} - V_{n_k}$,
satisfies the requirements.
\end{proof}

Our final two lemmas concern the properties of the limit random
distribution of determining sequences.

\begin{lemma} \label{lemma6}
Let $(X_n)_{n \geq 1}$ be a determining sequence of random variables with limit
random distribution $F$. Then for any set $A \subset \Omega$
with $P(A) > 0$ we have
\begin{equation}\label{2.7}
{\mathbb E}_A \left( \int\limits^{+\infty}_{-\infty} x^2 dF_\bullet (x) \right) =
\int\limits^{+\infty}_{-\infty} x^2 dF_A(x),
\end{equation}
in the sense that if one side is finite then the other side is
also finite and the two sides are equal. The statement remains
valid if in \eqref{2.7} we replace the interval of integrations by
$(-t, t)$, provided $t$ and $-t$ are continuity points of $F_A$.
\end{lemma}

\begin{proof}This lemma follows easily from \eqref{1.4} by integration by parts. \end{proof}

\begin{lemma} \label{lemma7}
Let $X, X_1, X_2, \dots$ be random variables such that both sequences
$(X_n)_{n \geq 1}$ and $(X_n - X)_{n \geq 1}$ are determining; let $F_\bullet$ and
$G_\bullet$ denote, respectively, their limit random
distributions. Then $\int\limits^{+\infty}_{-\infty} x^2 dG_\bullet (x)
< + \infty$ a.s.\ implies $\int\limits^{+\infty}_{-\infty} x^2
dF_\bullet (x) < +\infty$ a.s.\ and conversely.
\end{lemma}

\begin{proof}
Let $\ve > 0$ and choose a set $A \subset \Omega$ such that
$P(A) \geq 1 - \ve$ and on $A$ both $X$ and
$\int\limits^{+\infty}_{-\infty} x^2 dF_\bullet (x)$ are bounded. Let
$F_A$ and $G_A$ denote the limit random distribution of $(X_n)_{n \geq 1}$
resp.\ $(X_n - X)_{n \geq 1}$ relative to $A$. Replacing $X_n$ by $X_n +
\tau_n$ where $\tau_n \to 0$ a.s.\ clearly does not change the
limit distributions $F_A, G_A, F_\bullet, G_\bullet$, and thus by
passing to a subsequence and using Lemma~\ref{lemma5} we can assume,
without loss of generality, that the $X_n$ are identically
distributed on $A$. Then
$$
{\mathbb E}_A X^2_1 = {\mathbb E}_A X^2_2 = \cdots = \int\limits^{+\infty}_{-\infty} x^2
dF_A (x),
$$
where the last integral is finite by the boundedness of
$\int\limits^{+\infty}_{-\infty} x^2 dF_\bullet (x)$ on $A$ and
Lemma~\ref{lemma6}. By Minkowski's inequality and the boundedness of $X$ on
$A$ it follows that ${\mathbb E}_A((X_n - X)^2)$ is also bounded, and thus
Fatou's lemma implies that
$$
\int\limits^{+\infty}_{-\infty} x^2 dG_A(x) \leq \liminf_{n \to \infty}
{\mathbb E}_A \left((X_n - X)^2 \right) < + \infty.
$$
Using Lemma~\ref{lemma6} again it follows that $\int\limits^{+\infty}_{-\infty}
x^2 dG_\bullet (x) < + \infty$ a.s.\ on $A$. As the measure of
$A$ can be chosen arbitrarily close to $1$, we get $\int\limits x^2
dG_\bullet (x) < + \infty$ a.s., as required.
\end{proof}

\subsection{Proof of the theorems}

We begin with the proof of Theorem~1. Using diagonalization and
Chebyshev's inequality it follows that if a sequence $(X_n)_{n \geq 1}$
satisfies (B), then it has a subsequence bounded in probability
and thus also a determining subsequence. By Lemma~\ref{lemma2} the same
conclusion holds if $(X_n)_{n \geq 1}$ satisfies (A). Thus to prove our
theorem it suffices to prove the equivalence of (A), (B), (C),
(D) for determining sequences $(X_n)_{n \geq 1}$. In what follows we shall
prove the implications
(A)$\implies$(C)$\implies$(D)$\implies$(B); since
(B)$\implies$(A) follows easily from Theorem~D in the previous section by diagonalization, this
will prove Theorem~1.\\

Assume that $(X_n)_{n \geq 1}$ is a determining sequence satisfying (A),
i.e.\ there exists a subsequence $(X_{n_k})_{k \geq 1}$ and random variables $X, Y$
with $Y > 0$ such that for any further subsequence $(X_{m_k})_{k \geq 1}$
of $(X_{n_k})_{k \geq 1}$ and any set $A \subset \Omega$ with $P(A) > 0$ we
have
\begin{equation}\label{3.1}
P_A \left( \frac{\sum\limits_{k=1}^N (X_{m_k} - X)}{Y \sqrt
N} < t \right) \to \Phi(t) \quad \text{for all } t.
\end{equation}
We claim that $(X_n)_{n \geq 1}$ satisfies (C). Clearly we can assume
without loss of generality that $(X_{n_k})_{k \geq 1} = (X_k)_{k \geq 1}$ and since
$(X_n - X)_{k \geq 1}$ contains a determining subsequence, we can assume
also that $(X_n - X)_{n \geq 1}$ itself is determining. Moreover, since
$(X_n - X)_{n \geq 1}$ satisfies (C) if and only if $(X_n)_{n \geq 1}$ does (see Lemma~\ref{lemma7}), we can
assume that $X = 0$. Assume indirectly that $(X_n)_{n \geq 1}$ does not satisfy
(C), i.e.\ there exists a set $B \subset \Omega$ with $P(B) > 0$
such that
\begin{equation} \label{3.2}
\lim_{t \to \infty} \int\limits^t_{-t} x^2 dF_\bullet(x) = +\infty
\quad \text{on } B.
\end{equation}
Then there exists a set $B^* \subset B$ with $P(B^*) \geq
P(B)/2$ such that on $B^*$ the random variable $Y$ is bounded and \eqref{3.1}
holds uniformly, i.e.\ there exists a constant $K > 0$ and a
numerical sequence $K_t \to +\infty$ such that
$$
\int\limits^t_{-t} x^2 dF_\bullet (x) \geq K_t \ \text{ and } \ Y \leq
K \ \text{ on } \ B^* .
$$
Also, $1 - F_\bullet(t) + F_\bullet(-t) \to 0$ a.s.\  as $t \to
\infty$, and thus we can choose a set $B^{**} \subset B^*$ with
$P(B^{**}) \geq P(B^*) / 2$ such that on $B^{**}$ the last
convergence relation holds uniformly, i.e.\ there exists a
positive numerical sequence $\ve_t \searrow 0$ such that
\begin{equation} \label{3.4}
1 - F_\bullet (t) + F_\bullet(-t) \leq \ve_t \ \text{ on } \
B^{**}.
\end{equation}
We show that there exists a subsequence $(X_{m_k})_{k \geq 1}$ of
$(X_{n_k})_{k \geq 1}$ such that \eqref{3.1} fails for $A = B^{**}$. Since our
argument will involve the sequence $(X_n)_{n \geq 1}$ only on the set
$B^{**}$, in the sequel we can assume, without loss of
generality, that $B^{**} = \Omega$. That is, we may assume that
\eqref{3.4} holds on the whole probability space.\\

Let $C$ be an arbitrary set in the probability space with $P(C)
> 0$. Integrating \eqref{3.4} and using \eqref{1.4} and Lemma~\ref{lemma6} we get
\begin{equation}\label{3.5}
\int\limits^t_{-t} x^2 dF_C(x) \geq K_t, \quad
1 - F_C(t) + F_C(-t) \leq\ve_t, \quad t \in H,
\end{equation}
where $H$ denotes the set of continuity points of $F_\Omega$.
Choose $t_0 \in H$ so large that $\ve_{t_0} \leq 1/16$ and then
choose $t_1$ so large that
$$
K^{1/2}_t /4  \geq 2t^2_0 \ \text{ for } \ t \geq t_1, \ t \in H.
$$
Then for $t \geq t_1$, $t \in H$ we have, using the second
relation of \eqref{3.5},
\begin{eqnarray*}
\left| \int\limits^t_{-t} x dF_C(x) \right| & \leq & 2t^2_0 + \int\limits_{t_0
\leq |x| \leq t} |x| dF_C(x) \\
&\leq & 2t^2_0 + \left( \int\limits_{|x| \geq t_0} dF_C(x) \right)^{1/2}
\left( \int\limits_{|x| \leq t} x^2 dF_C(x) \right)^{1/2} \\
&\leq & 2t^2_0 + \frac{1}{4} \left( \int\limits_{|x| \leq t} x^2 dF_C(x)
\right)^{1/2} \\
& \leq & \frac{1}{2} \left(\int\limits_{|x| \leq t} x^2
dF_C(x) \right)^{1/2},
\end{eqnarray*}
and thus we proved that for any $C \subset \Omega$ with $P(c) >
0$ we have
\begin{equation}\label{3.6}
\int\limits^t_{-t} x^2 dF_C(x) - 2 \left( \int\limits^t_{-t}
xdF_C(x)\right)^2 \geq \frac{1}{2} K_t,
\qquad t \geq t_1, \ t \in H.
\end{equation}
Since $(X_n)_{n \geq 1}$ is bounded in probability, there exists a function
$\psi(x) \nearrow \infty$ such that
\begin{equation}\label{3.7}
\sup_n E \psi(X_n) \leq 1
\end{equation}
(see [9]). Let $(a_k)_{k \geq 1}$ be a sequence of integers tending to
$+\infty$ so slowly that $a_k \leq \log k$ and
\begin{equation}\label{3.8}
\delta_k : = a_k / \psi(k^{1/4}) \to 0.
\end{equation}
Let further $(\ve_n)_{n \geq 1}$ tend to $0$ so rapidly that
$\ve_{a_k} \leq 2^{-k}$. By Theorem~\ref{th8d} there exists a
subsequence $(X_{m_k})_{k \geq 1}$ and a sequence $(Y_k)_{k \geq 1}$ of discrete
random variables such that \eqref{2.1B} holds and for each $k > 1$ the atoms of
the finite $\sigma$-field $\sigma \{ Y_1, \dots, Y_{a_k} \}$ can
be divided into two classes $\Gamma_1$ and $\Gamma_2$ such that
\begin{equation}\label{3.9}
\sum_{B \in \Gamma_1} P(B) \leq \ve_{a_k} \leq 2^{-k}
\end{equation}
and for each $B \in \Gamma_2$ there exist i.i.d.\ random variables
$Z^{(B)}_{a_k + 1}, \dots, Z^{(B)}_k$ on $B$ with distribution
$F_B$ such that
\begin{equation}\label{3.10}
P_B \bigl( |Y_j - Z^{(B)}_j | \geq 2^{-k} \bigr) \leq 2^{-k},
\qquad j = a_k + 1, \dots, k.
\end{equation}
We show that
\begin{equation}\label{3.11}
P \left( \frac{\sum\limits_{k=1}^N X_{m_k}}{Y \sqrt N} <
t\right) \to \Phi(t) \quad \text{for all } t
\end{equation}
cannot hold; this will complete our indirect proof of
(A)$\implies$(C). Set
\begin{align*}
S^{(B)}_{a_k, k} &= \sum^k_{j = a_k + 1} Z^{(B)}_j, \qquad B \in
\Gamma_2, \\
\overline S_{a_k, k} &= \sum_{B \in \Gamma_2} S^{(B)}_{a_k, k} \mathbbm{1}_B,
\end{align*}
where $\mathbbm{1}_B$ denotes the indicator function of $B$.
By \eqref{3.10},
$$
P_B \left( \left| \sum^k_{j = a_k + 1} Y_j - \sum^k_{j = a_k +
1} Z^{(B)}_j \right| \geq 1 \right) \leq 2^{-a_k}, \qquad B \in
\Gamma_2,
$$
and thus using \eqref{3.9} we get
\begin{equation}\label{3.12}
P \left( \left| \sum^k_{j = a_k + 1} Y_j - \overline S_{a_k,
k} \right| \geq 1 \right) \leq 2 \cdot 2^{-k}.
\end{equation}
By \eqref{3.7}, \eqref{3.8} and the Markov inequality we have
\begin{eqnarray*}
P \left( \left| \sum^{a_k}_{j = 1} X_{m_j} \right| \geq a_k
k^{1/4} \right) & \leq & a_k \sup_{1 \leq j \leq a_k} P \left(
|x_{m_j}| \geq k^{1/4} \right) \\
& \leq & a_k \psi(k^{1/4})^{-1} \\
& = & \delta_k,
\end{eqnarray*}
which, together with \eqref{3.12} and \eqref{2.1B}, yields
\begin{equation}\label{3.13}
P \left( \left| \sum^k_{j = 1} X_{m_j} - \overline S_{a_k, k}
\right| \geq 2_{a_k} k^{1/4} \right) \leq 3 \cdot 2^{-a_k} + \delta_k.
\end{equation}
Applying Lemma~\ref{lemma4} to the i.i.d.\ sequence $\{ Z^{(B)}_j, a_k + 1
\leq j \leq k \}$ and using \eqref{3.6} we get
\begin{align*}
P_B \left(\left| \frac{S^{(B)}_{a_k, k}}{\sqrt k} \right|
\leq 1 \right)
&\leq P_B \left( \frac{S^{(B)}_{a_k, k}}{\sqrt{k -
a_k}} \leq 2 \right)  \\
&\leq \text{\rm const.}\, K^{-1/2}_{\sqrt k},
\end{align*}
where the constant is absolute. Thus using \eqref{3.9} and $Y \leq K$
it follows that
\begin{equation}\label{3.14}
P \left( \left| \frac{\overline S_{a_k, k}}{Y \sqrt k} \right|
\leq \frac{1}{K} \right) \leq \text{\rm const.} \,
K^{-1/2}_{\sqrt k} + 2^{-k}.
\end{equation}
If \eqref{3.11} were true then by \eqref{3.13} and $a_k \leq \log k$ we
would also have
$$
\overline S_{k, a_k} / Y \sqrt k \overset {\mathcal D}
\longrightarrow N(0,1),
$$
which clearly contradicts \eqref{3.14} for large $k$, since the
right-hand side tends to zero. This completes the proof of
(A)$\implies$(C).\\

The remaining implications (C)$\implies$(D) and (D)$\implies$(B)
of Theorem~\ref{th1} are easy. Assume first that (C) holds, then for any
$\ve > 0$ there exists a set $A \subset \Omega$ with $P(A) \geq
1 - \ve$ and a constant $K = K_\ve$ such that
$$
\int\limits^{+ \infty}_{-\infty} x^2 dF_\bullet (x) \leq K \ \text{ on
}\ A.
$$
Integrating the last relation on $A$ and using Lemma~\ref{lemma6} we get
\eqref{1.5}, i.e.\ (D) holds. Assume now that (D) holds, i.e.\ for any
$\ve > 0$ there exists a set $A \subset \Omega$ with $P(A) \geq
1 - \ve$ such that \eqref{1.5} is valid. Applying Lemma~\ref{lemma5} for $(X_n)_{n \geq 1}$
on the set $A$ it follows that there exists a subsequence
$(X_{n_k})_{k \geq 1}$ and random variables $Y_k$ and $\tau_k,~ k = 1,2,\dots$,
defined on $A$ such that $X_{n_k} = Y_k + \tau_k, ~k =
1,2,\dots$, and such that the random variables $Y_k$ have distribution $F_A$ on $A$ and
$\tau_k \to 0$ a.s.\ on $A$. Choose a set $B \subset A$ with
$P(B) \geq 1 - 2\ve$ such that $\tau_k \to 0$ uniformly on~$B$.
Then clearly $(\tau_n)_{n \geq 1}$ is uniformly bounded on $B$, and further
\begin{eqnarray*}
\int\limits_B Y^2_k dP &\leq & \int\limits_A Y^2_k dP \\
& = & P(A)
\int\limits^{+\infty}_{-\infty} x^2 dF_A(x) \\
&\leq & \int\limits^{+\infty}_{-\infty} x^2 dF_A(x) < + \infty
\end{eqnarray*}
for each $k \geq 1$ by the identical distribution of the $Y_k$'s
and \eqref{1.5}. Thus on $B$ the sequences $(Y_k)_{k \geq 1}$ and $(\tau_k)_{k \geq 1}$ have
bounded $L_2$ norms and thus the same holds for $X_{m_k} = Y_k +
\tau_k$, i.e.\
$$
\sup_k \int\limits_B X^2_{m_k} dP < +\infty.
$$
In view of $P(B) \geq 1 - 2\ve$ this shows that $(X_n)_{n \geq 1}$
satisfies statement (B). This completes the proof of Theorem~1.\\

{\it Proof of Theorem 2}.
Theorem~\ref{th2} follows from Theorem~\ref{th1} and a slightly sharper
form of Theorems~\ref{thb} and~\ref{thc} which was proved in~\cite{kape}. We already mentioned the fact
that in Theorem~\ref{thb} the random variables $X,Y$ appearing in \eqref{1.2} actually
satisfy $X \in L_2$, $Y \in L_{1/2}$. Moreover, if instead of
\eqref{1.1} we make the slightly stronger assumption that the sequence
$(X^2_n)_{n \geq 1}$ is uniformly integrable then by the weak compactness
criteria in $L_1$ and $L_2$ it follows that there exists a
subsequence $(X_{n_k})$ and random variables $X \in L_2$, $Y \in L_{1/2}$
such that
\begin{equation}\label{3.15}
X_{n_k} \to X \text{ weakly in } L_2, \qquad
(X_{n_k} - X)^2 \to Y^2 \  \text{weakly in } L_1.
\end{equation}
As is shown in~\cite{lan}, in this case \eqref{1.2} holds with the random variables
$X, Y$ determined by \eqref{3.15}. We turn now to the proof of Theorem~2. As in the case of
Theorem~1, it suffices to prove the equivalence of statements
(A), (B), (C) in the case when $(X_n)_{n \geq 1}$ is determining. Also,
since replacing $X_n$ by $X_n + \tau_n$ where $\sum |\tau_n| < +
\infty$ a.s.\ does not affect the validity of \eqref{1.2}, the
conclusion (B)$\implies$(A) of Theorem \eqref{th2} is contained in the
stronger form of Theorem~\ref{thb} mentioned above. Thus it suffices to
verify the implications (A)$\implies$(C) and (C)$\implies$(B).
To prove (A)$\implies$(C) let us assume that $(X_n)_{n \geq 1}$ is
determining with limit distribution $F$, and that there exist a
subsequence $(X_{n_k})_{k \geq 1}$ and random variables $X \in L_2$, $Y \in L_{1/2}$,
$Y > 0$, such that for all subsequences $(X_{m_k})_{k \geq 1}$ of
$(X_{n_k})_{k \geq 1}$ and any set $A \subset \Omega$ with $P(A) > 0$ equation \eqref{1.2}
holds. We show $\int\limits^{+\infty}_{-\infty} x^2 dF(x) < + \infty$.
Clearly we can assume without loss of generality that $(X_{n_k})_{k \geq 1}
< (X_k)_{k \geq 1}$. Fix $\ve > 0$. By the implication (A)$\implies$(B)$\implies$(D)
of Theorem~\ref{th1} there is a set $A \subset\Omega$ with $P(A) \geq 1
- \ve$ and a subsequence $(X_{n_k})_{k \geq 1}$ such that
\begin{equation}\label{3.16}
\sup_k \int\limits_A X^2_{n_k} dP < + \infty \ \text{ and } \
\int\limits_A x^2 dF_A(x) < + \infty,
\end{equation}
where $F_A$ is the limit distribution of $(X_n)_{n \geq 1}$ on $A$.
Applying Lemma~\ref{lemma5} for $(X_{n_k})_{k \geq 1}$ on $A$ it follows that there
exists a subsequence $(X_{m_k})_{k \geq 1}$ of $(X_{n_k})_{k \geq 1}$ admitting the decomposition
\begin{equation}\label{3.17}
X_{m_k} = Y_k + \tau_k \ \text{ on } \ A,
\end{equation}
where the $Y_k$ are identically distributed on $A$ with
distribution function $F_A$ and $\sum |\tau_k| < + \infty$ a.s.\
on $A$. Being an identically distributed sequence with finite
expectation, the sequence $(Y^2_k)_{k \geq 1}$ is uniformly integrable on
$A$, and thus by the sharper form of Theorem~\ref{thb} mentioned above it
follows that there exists a subsequence $(Y_{p_k})_{k \geq 1}$ of $(Y_k)_{k \geq 1}$
such that
$$
Y_{p_k} \to U \text{ weakly in } L_2(A), \qquad
(Y_{p_k} - U)^2 \to V^2 \text{ weakly in } L_1(A),
$$
and for any $B \subset A$ with $P(B) > 0$ we have
$$
P_B \left( \frac{\sum_{k=1}^N (Y_{p_k} - U)}{V \sqrt
N} < t \right) \to \Phi (t) \ \text{ for all } t,
$$
where $U$, $V$ are random variables such that $U \in L_2(A)$, $V \in
L_{1/2}(A)$, $V > 0$. Thus by \eqref{3.17} and $\sum |\tau_k| < +
\infty$ a.s.\ on $A$ we get
$$
P_B \left( \frac{\sum_{k=1}^N (X_{m_{p_k}} - U)}{V \sqrt
N} < t \right) \to \Phi (t) \ \text{ for all } t
$$
for any $B \subset A$ with $P(B) > 0$. Comparing with \eqref{1.2} and
using Lemma~\ref{lemma3} we get $U = X$, $V = Y$ a.s.\ on $A$, and thus we proved that
$$
Y_{p_k} \to X \text{ weakly in } L_2(A), \qquad
(Y_{p_k} - X)^2 \to Y^2 \text{ weakly in } L_1(A).
$$
Hence
\begin{align}
{\mathbb E}_A Y^2 &= \lim_{k \to \infty} {\mathbb E}_A(Y_{p_k} - X)^2 = \lim_{k \to
\infty} ({\mathbb E}_A Y^2_{p_k} - 2 {\mathbb E}_A Y_{p_k} X + {\mathbb E}_A X^2) \nonumber \\
&= \lim_{k \to \infty} {\mathbb E}_A Y^2_{p_k} - {\mathbb E}_A X^2 =
\int\limits^{+\infty}_{-\infty} x^2 dF_A(x) - {\mathbb E}_A X^2, \label{3.18}
\end{align}
where in the last step we used the fact that the $Y_k$'s have
distribution $F_A$ on $A$. Hence
\begin{equation}\label{3.19}
P(A)^{-1} \int\limits_A Y^2 dP = \int\limits^{+\infty}_{-\infty} x^2 dF_A(x) -
P(A)^{-1} \int\limits_A X^2 dP .
\end{equation}
Since $X \in L_2(\Omega)$, $Y^2 \in L_1(\Omega)$, the left-hand
side of \eqref{3.19} and the second term on the right-hand side
approach finite limits as $P(A) \to 1$ and thus
$\int^{+\infty}_{-\infty} x^2 dF_A(x)$ also converges to a
finite limit. On the other hand, $F_A \to F$ as $P(A) \to 1$ and
thus Fatou's lemma implies
$$
\int\limits^{+\infty}_{-\infty} x^2 dF(x) \leq \liminf_{P(A) \to 1}
\int\limits^{+\infty}_{-\infty} x^2 dF_A(x) < + \infty,
$$
proving the implication (A)$\implies$(C). Now if (C) holds then
by Lemma~\ref{lemma5} there exists a subsequence $(X_{n_k})_{k \geq 1}$ permitting the
decomposition \eqref{1.7} where $\sum |\tau_k| < +\infty$ a.s.\ and
$Y_k$ are identically distributed with distribution $F$; since
$F$ has finite variance by (C), the first relation of \eqref{1.8}
holds. Thus $(X_n)_{n \geq 1}$ satisfies (B) and the proof of Theorem~\ref{th2} is
completed.

\section*{Acknowledgments}

Christoph Aistleitner is supported by the Austrian Science Fund (FWF), projects F-5512, I-3466, I-4945, I-5554, P-34763, P-35322 and Y-901. Istvan Berkes is supported by Hungarian Foundation NKFI-EPR No.\ K-125569.

\end{document}